\documentclass[12pt,twoside]{article}
\usepackage{amssymb}
\renewcommand{\leq}{\leqslant}
\renewcommand{\geq}{\geqslant}

\begin{document}
\parindent=0.3in
\parskip=0in
\baselineskip=20pt plus 1pt

\vspace*{0.1in}
\begin{center}
{\large 
THE DECOMPOSITION OF LIE POWERS}
\end{center}

\vspace{0.01in}
\begin{center}
R. M. BRYANT {\it and\/} M. SCHOCKER
\end{center}

\footnote{\ Research of second author supported by Deutsche
Forschungsgemeinschaft (DFG Scho-799).}

\begin{center}
{\sc Abstract}
\end{center}

\noindent
{\small Let $G$ be a group, $F$ a field of prime
characteristic $p$ and $V$ a finite-dimensional $FG$-module.
Let $L(V)$ denote the free Lie algebra on $V$ regarded as
an $FG$-submodule of the free associative algebra (or tensor
algebra) $T(V)$. For each positive integer $r$, let $L^r(V)$
and $T^r(V)$ be the $r$th homogeneous components of $L(V)$
and $T(V)$, respectively. Here $L^r(V)$ is called the $r$th
Lie power of $V$. Our main result is that there are submodules
$B_1$, $B_2$, \dots\ of $L(V)$ such that, for all $r$, $B_r$
is a direct summand of $T^r(V)$ and, whenever $m \geq 0$ and
$k$ is not divisible by $p$,
$$L^{p^mk}(V) = L^{p^m}(B_k) \oplus L^{p^{m-1}}(B_{pk})
\oplus \cdots \oplus L^p(B_{p^{m-1}k}) \oplus L^1(B_{p^mk}).$$

}
\noindent
{\small Thus every Lie power is a direct sum of Lie powers of
$p$-power degree. The approach builds on an analysis of $T^r(V)$
as a bimodule for $G$ and the Solomon descent algebra.

}
\vspace{0.4in}
\baselineskip=20pt plus 1pt
\begin{center}
1. {\it Introduction}
\end{center}
\renewcommand{\theequation}{1.\arabic{equation}}
\setcounter{equation}{0}

\noindent
Let $G$ be a group and $F$ a field. For any finite-dimensional
$FG$-module $V$ let $L(V)$ be the free Lie algebra on $V$. Then
$L(V)$ may be regarded as an $FG$-module on which each element
of $G$ acts as a Lie algebra automorphism. Furthermore, each
homogeneous component $L^r(V)$ is a finite-dimensional submodule,
called the $r$th Lie power of $V$.  We regard $L(V)$ as an
$FG$-submodule of the free associative algebra, or tensor algebra,
$T(V)$. Thus $L^r(V)$ is a submodule of $T^r(V)$.

In characteristic $0$, the structure of $L^r(V)$ has been clarified
by a number of papers, including those of Brandt [{\bf 3}],
Klyachko [{\bf 19}], and Kraskiewicz and Weyman [{\bf 20}]. 
We concentrate here on the more
difficult case where $F$ has prime characteristic $p$. When
$r$ is not divisible by $p$, $L^r(V)$ is comparatively well
understood. In this case, $L^r(V)$ is a direct summand of
$T^r(V)$ (that is, $L^r(V)$ has a complement in $T^r(V)$)
and it is possible to exploit knowledge of $T^r(V)$.
For example, when $F$ is
infinite and $G$ is the general linear group on $V$, the
indecomposable summands of $L^r(V)$, for $r$ not divisible
by $p$, were described up to isomorphism by Donkin and
Erdmann [{\bf 13}].

The first case beyond this, namely the case $r=p$, was
studied by Bryant and St\"ohr [{\bf 10}], by means of a detailed
analysis of $T^p(V)$.
Further progress was made by Erdmann and Schocker [{\bf 14}] who
studied the case $r=pk$
where $k$ is not divisible by $p$. Their methods make substantial
use of Solomon's descent algebra, showing that its significance
for free Lie algebras goes beyond the case of characteristic $0$.
Their main results are concerned with the relationship
between $L^{pk}(V)$ and $L^k(V)$. As a consequence, they were
able to prove the factorisation $\Phi_{FG}^{pk} =
\Phi_{FG}^p \circ \Phi_{FG}^k$ for the ``Lie resolvents" introduced
by Bryant [{\bf 4}].
Such factorisations, for certain groups $G$, first appeared in
[{\bf 5}] and [{\bf 6}]. It follows from the results of [{\bf 14}] that the
study of $L^{pk}(V)$ can be reduced, broadly speaking, to the study
of $L^p(L^k(V))$.  As we have already seen, $L^k(V)$ is
comparatively well understood. Thus the study of $pk$th Lie
powers is reduced, in a sense, to the study of $p$th Lie powers.

This theme is developed a great deal further in the present paper.
We show that the study of arbitrary Lie
powers may be reduced, in the same sense, to the study of Lie powers
of $p$-power degree.

First we develop some aspects of the Poincar\'e--Birkhoff--Witt
Theorem. Linked with this theorem there is a filtration of $T^r(V)$
by $FG$-submodules in which the factors are isomorphic to
certain tensor products of symmetric powers of Lie powers of $V$,
called ``higher Lie modules".  At the bottom of the filtration is
$L^r(V)$ and, at the top, the $r$th symmetric power 
$S^r(V)$. Several papers on free Lie algebras ([{\bf 10}], for example)
have used this filtration to study $L^r(V)$, exploiting
the fact that $T^r(V)$ is comparatively well understood.
The weakness of this method
is that the filtration has many terms.

Our first result, called the ``Filtration Theorem" and
proved in \S3, does something to remedy this weakness. We
show that the filtration can be modified by means of primitive
idempotents of the modular descent algebra. It turns out that
$T^r(V)$ splits up as a direct sum of $FG$-submodules each with
a much shorter filtration than before. The summands of
$T^r(V)$ correspond to the primitive idempotents, and these are
indexed by the ``$p$-equivalence classes" of partitions of $r$.

In \S4, we prove our main result about free Lie algebras, called
the ``Decomposition Theorem". This reduces the study
of arbitrary Lie powers to the study of Lie powers of
$p$-power degree and may be stated as follows. Let $F$ be a
field of characteristic $p$. Let $G$ be a group and let $V$ be
a finite-dimensional $FG$-module. Then, for each positive
integer $r$, there is a submodule $B_r$ of $L^r(V)$ such that
$B_r$ is a direct summand of $T^r(V)$ and, for $k$ not
divisible by $p$ and $m \geq 0$,
$$L^{p^mk}(V) = L^{p^m}(B_k) \oplus L^{p^{m-1}}(B_{pk})
\oplus \cdots \oplus L^1(B_{p^mk}).$$
This has to be interpreted with a little care. As we shall see in \S2,
the subalgebra of $L(V)$ generated by $B_r$ is free on $B_r$. It may
therefore be identified with $L(B_r)$ and then we have
$L^i(B_r) \subseteq L^{ir}(V)$ for all $i$. The proof of the
Decomposition Theorem uses the Filtration Theorem of \S3
as well as Lazard's ``Elimination Theorem" (see \S2 below).
This last result has proved to be
fundamental in the theory and has been repeatedly used in
earlier work such as [{\bf 7}].

In another paper [{\bf 9}] we shall obtain further information about the
isomorphism types of the modules $B_r$ and we shall show that the
Decomposition Theorem yields the factorisation of Lie resolvents
$\Phi_{FG}^{p^mk} = \Phi_{FG}^{p^m} \circ \Phi_{FG}^k$
for every non-negative integer $m$ and every positive integer
$k$ not divisible by $p$.

\vspace{0.4in}
\begin{center}
2. {\it Preliminaries}
\end{center}
\renewcommand{\theequation}{2.\arabic{equation}}
\setcounter{equation}{0}

We shall require some technical variations on known material concerned
with free Lie algebras and representations of the
general linear group. These are collected together in the
present section.

We first work over a ring $R$ which is allowed to be either the
ring of integers or an arbitrary field. (This is simply a sufficient
degree of generality for our purposes.)
Let $V$ be a free $R$-module.
We write $T(V)$ for the free associative $R$-algebra
(with identity element) freely generated by any $R$-basis of $V$.
Thus $T(V)$ has an $R$-module decomposition
$T(V) = \bigoplus_{r \geq 0} T^r(V)$ 
where, for each $r$, $T^r(V)$ is the $r$th homogeneous component
of $T(V)$ and we identify $T^1(V)$ with $V$. Each $T^r(V)$ is
a free $R$-module and
$T^r(V) \cong V \otimes_R \cdots \otimes_R V$, with $r$ factors.
We call $T^r(V)$ the $r$th tensor power of $V$.
Note, however, that we write $ab$ rather than $a \otimes b$ for
the product of elements $a$ and $b$ of $T(V)$. For $a \in T^r(V)$
we write $\deg a = r$.

The algebra $T(V)$ may be made into a Lie $R$-algebra by means of the
bracket operation given by $[a,b] = ab - ba$. The Lie subalgebra of
$T(V)$ generated by $V$ is denoted by $L(V)$. By a theorem
of Witt (see [{\bf 18}] for the case of a field and [{\bf 2}] for general
$R$), $L(V)$ is a free Lie $R$-algebra, freely
generated by any basis of $V$, and $T(V)$ is the enveloping algebra
of $L(V)$. For $r \geq 1$, we write $L^r(V) = T^r(V) \cap L(V)$.
Thus $L(V) = \bigoplus_{r \geq 1} L^r(V)$, where $L^1(V) = V$.
Each $L^r(V)$ is a free $R$-module, called the $r$th Lie power
of $V$.

Let $S(V)$ denote the free commutative associative $R$-algebra
(that is, polynomial ring) freely generated by any $R$-basis of $V$.
Then, as for $T(V)$, we have $S(V) = \bigoplus_{r \geq 0}S^r(V)$.
Each $S^r(V)$ is a free $R$-module called the $r$th symmetric power
of $V$.

When $V$ has finite rank then $T^r(V)$, $L^r(V)$ and
$S^r(V)$ also have finite rank.

Suppose that $G$ is a group and that the free $R$-module $V$ is a
(right) $RG$-module.
The action of $G$ on $V$ extends to $T(V)$ so that $G$ acts by algebra
automorphisms. Thus $T(V)$ becomes an $RG$-module. The subspaces
$L(V)$, $T^r(V)$ and $L^r(V)$ are $G$-invariant, so they also become
$RG$-modules. Similarly, $S(V)$ and the subspaces $S^r(V)$ can be
given the structure of $RG$-modules.

We now consider the case $R=F$, where $F$ is a field.
Suppose that $B$ is a subspace of $T^r(V)$ for some $r \geq 1$. Then
it is easily verified that the subalgebra of $T(V)$ generated by $B$ is
isomorphic to $T(B)$ and may be identified with $T(B)$ by means
of an isomorphism which is the identity on $B$. Then 
$T^i(B) \subseteq T^{ir}(V)$ for all $i \geq 0$. Furthermore,
if $B \subseteq L^r(V)$ then the Lie subalgebra of $L(V)$ generated
by $B$ is equal to the subalgebra $L(B)$ of $T(B)$ as previously
defined. Thus it is a free Lie algebra. Furthermore, 
$L^i(B) \subseteq L^{ir}(V)$
for all $i \geq 1$. This is the justification for the notation
used in the Decomposition Theorem as stated in \S1.

If $L$ is a free Lie algebra over $F$ and
$W$ is a subspace of $L$ such that some basis (equivalently,
every basis) of $W$ is a free generating set of $L$ then we
say (with slight abuse of language) that $L$ is free on $W$
or that $L$ is freely generated by $W$. In this case $L$ is
isomorphic to $L(W)$.

If $L$ is a free Lie algebra over $F$ and $Q$ is a subalgebra
of $L$ then, by the theorem of Shirshov and Witt (see [{\bf 21},
Theorem 2.5]), $Q$ is free (on some free generating set). If
$W$ is a subspace of $L$ then the subalgebra of $L$
generated by $W$ is not necessarily free on $W$ (though it
is free on some subspace). If this subalgebra is free
on $W$ we write it as $L(W)$. Indeed we identify the subalgebra
generated by $W$ with $L(W)$ by means of the isomorphism which
is the identity on $W$.

The following result is essentially well known and is given
as [{\bf 8}, Lemma 1].

\vspace{0.2in}
{\sc Lemma} 2.1.  {\it Let\/ $L$ be a free Lie algebra over\/
$F$ and, for each\/ $r \geq 1$, let\/ $L_r$ be the\/ $r$th
homogeneous component of\/ $L$ with respect to a given free
generating set.  Let\/ $Q$ be a subalgebra of\/ $L$ of the
form\/ $Q = Q_1 \oplus Q_2 \oplus \cdots\,$ with\/
$Q_r = Q \cap L_r$ for all\/ $r$. For each\/ $r \geq 1$, let\/
$Q(r)$ denote the subalgebra of\/ $Q$ generated by\/
$Q_1 \oplus \cdots \oplus Q_r$, and write\/ $Q(0) = 0$.
For each\/ $r \geq 1$, let\/ $W_r$ be any subspace of\/ $Q_r$
satisfying\/ $Q_r = (Q(r-1) \cap Q_r) \oplus W_r$. Then\/ $Q$
is free on\/ $W_1 \oplus W_2 \oplus \cdots \;$, that is,
$$Q = L(W_1 \oplus W_2 \oplus \cdots).$$}

\vspace{0.01in}
For subspaces $B$ and $C$ of any Lie algebra over $F$, let $C \wr B$
denote the subspace defined by
$$C \wr B = C + [C,B] + [C,B,B] + \cdots\; .$$
(We use the left-normed convention for Lie products.)
The following lemma is a version of Lazard's ``Elimination
Theorem" (see [{\bf 2}, Chapter 2, \S2.9, Proposition 10]). In the
form written here it is a special case (with trivial group action)
of [{\bf 7}, Lemma 2.2] or [{\bf 8}, Lemma 2].

\vspace{0.2in}
{\sc Lemma} 2.2. {\it Let\/ $B$ and $C$ be\/
$F$-spaces, and consider the
free Lie algebra\/ $L(B \oplus C)$. Then\/ $B$ and\/ $C \wr B$
freely generate subalgebras\/ $L(B)$ and\/ $L(C \wr B)$, and
there is an\/ $F$-space decomposition
$L(B \oplus C) = L(B) \oplus L(C \wr B)$. Furthermore,
$$C \wr B = C \oplus [C,B] \oplus [C,B,B] \oplus \cdots$$
and, for each\/ $n \geq 0$, there is an\/ $F$-space isomorphism
$$\zeta_n: [C,\underbrace{B,\dots,B}_n] \stackrel{\cong}{\longrightarrow}
C \otimes \underbrace{B \otimes \cdots \otimes B}_n$$
such that\/ $\zeta_n([c,b_1,\dots,b_n]) =
c \otimes b_1 \otimes \cdots \otimes b_n$ for all\/ $c \in C$
and\/ $b_1,\dots,b_n \in B$.}

\vspace{0.2in}
{\sc Corollary} 2.3.  {\it Suppose, in Lemma\/ {\rm 2.2}, that\/
$C$ is the direct sum of subspaces, $C = \bigoplus_i X_i$.
Then\/  $C \wr B = \bigoplus_i (X_i \wr B)$ and, for each\/ $i$,
$$X_i \wr B = X_i \oplus [X_i,B] \oplus [X_i,B,B] \oplus \cdots\; .$$
Furthermore, for each\/ $n$, $\zeta_n$ restricts to an isomorphism
$$[X_i,\underbrace{B,\dots,B}_n] \stackrel{\cong}{\longrightarrow}
X_i \otimes \underbrace{B \otimes \cdots \otimes B}_n.$$}

\vspace{0.01in}
We next consider Schur algebras for the general linear group.
Basically, we follow [{\bf 16}], but we adapt the treatment there
to allow the underlying field to be finite or infinite. Details
omitted in our treatment can be found in [{\bf 16}] or are
straightforward modifications of what can be found there.
As before, $R$ is either the ring of integers or an arbitrary field.

Let $n$ and $r$ be positive integers. Let $I(n,r)$ be the set of
all ordered $r$-tuples ${\bf i} = (i_1,\dots,i_r)$ where
$i_1,\dots,i_r \in \{1,\dots,n\}$.

Let $A_R(n,r)$ be the homogeneous component of degree $r$ in the
polynomial ring over $R$ in $n^2$ indeterminates
$c_{ij}$ ($1 \leq i,j \leq n$). Thus $A_R(n,r)$ has an $R$-basis consisting
of the monomials of degree $r$.  For ${\bf i}, {\bf j} \in I(n,r)$
where ${\bf i} = (i_1,\dots,i_r)$, ${\bf j} = (j_1,\dots,j_r)$,
we write $c_{{\bf i},\, {\bf j}}$ for the monomial
$c_{i_1j_1}\cdots c_{i_rj_r}$. The $c_{{\bf i},\, {\bf j}}$ are not
distinct (when $n,r > 1$) but they give a basis (with repetitions)
of $A_R(n,r)$.

Let $S_R(n,r) = {\rm Hom}_R(A_R(n,r),R)$. Thus $S_R(n,r)$ has an
$R$-basis (with repetitions) consisting of the elements
$\xi_{{\bf i},\,{\bf j}}$ (${\bf i}, {\bf j} \in I(n,r)$) where
$\xi_{{\bf i},\,{\bf j}}(c_{{\bf i},\,{\bf j}}) =1$ and
$\xi_{{\bf i},\,{\bf j}}(c_{{\bf i}',\,{\bf j}'}) = 0$ if
$c_{{\bf i}',\,{\bf j}'} \neq c_{{\bf i},\,{\bf j}}$.  Multiplication
in $S_R(n,r)$ may be defined as in [{\bf 16}, \S2.3]: for
$\xi, \eta \in S_R(n,r)$,
$$(\xi\eta)(c_{{\bf i},\,{\bf j}}) = \sum_{{\bf k} \in I(n,r)}
\xi(c_{{\bf i},{\bf k}}) \eta(c_{{\bf k},\,{\bf j}}).$$
Thus $S_R(n,r)$ becomes an associative $R$-algebra with identity element.
This is the Schur algebra of degree $r$. As in [{\bf 16}, (2.3b)]
there are equations
$$\xi_{{\bf i},\,{\bf j}} \xi_{{\bf k},{\bf l}} =
\sum z_{{\bf s},{\bf t}} \xi_{{\bf s},{\bf t}},$$
where the $z_{{\bf s},{\bf t}}$ are non-negative integers depending
on ${\bf i}, {\bf j}, {\bf k}, {\bf l}$ but independent of $R$. When
we need to show the role of $R$ we write $\xi_{{\bf i},\,{\bf j}}^R$
instead of $\xi_{{\bf i},\,{\bf j}}$. We may identify $S_R(n,r)$ with
$R \otimes_{\mathbb Z} S_{\mathbb Z}(n,r)$ by identifying
$\xi_{{\bf i},\,{\bf j}}^R$ with $1 \otimes
\xi_{{\bf i},\,{\bf j}}^{\mathbb Z}$ for all ${\bf i},{\bf j}$.
If $E$ is a field containing $R$ we may regard $S_R(n,r)$ as a
subring of $S_E(n,r)$ by identifying $\xi_{{\bf i},\,{\bf j}}^R$ with
$\xi_{{\bf i},\,{\bf j}}^E$. We may also
identify $S_E(n,r)$ with $E \otimes_R S_R(n,r)$.

For $w \in A_R(n,r)$ and $g = (g_{ij}) \in {\rm GL}(n,R)$ we write
$w(g)$ for the element of $R$ obtained from $w$ by evaluating $c_{ij}$
as $g_{ij}$ for all $i,j$. We write $e(g)$ for the element of
$S_R(n,r)$ given by $e(g)(w) = w(g)$ for all $w \in A_R(n,r)$. Then,
as in [{\bf 16}, \S2.4], $e$ extends to an algebra homomorphism
\begin{equation}
e : R{\rm GL}(n,R) \longrightarrow S_R(n,r).
\end{equation}
This is surjective if $R$ is an infinite field, by [{\bf 16}, (2.4b)].

We shall consider right modules for $R{\rm GL}(n,R)$ and
$S_R(n,r)$, whereas [{\bf 16}] has left modules. This leads only
to minor, mainly notational, differences. Suppose that $U$ is an
$S_R(n,r)$-module which is free of finite rank as an $R$-module.
Then $U$ may be regarded as an
$R{\rm GL}(n,R)$-module by means of $e$, and such $R{\rm GL}(n,R)$-modules
are called polynomial modules of degree $r$. If $U$ has
$R$-basis $\{u_1,\dots,u_d\}$ then there are uniquely
determined elements $w_{st}$ ($1 \leq s,t \leq d$) of $A_R(n,r)$
such that, for all $\xi \in S_R(n,r)$ and all $s$,
\begin{equation}
u_s\xi = \sum_t \xi(w_{st})u_t.
\end{equation}
These polynomials $w_{st}$ are called the coefficient polynomials
(with respect to the given basis). For all $g \in {\rm GL}(n,R)$
we also have
\begin{equation}
u_sg = \sum_t w_{st}(g) u_t.
\end{equation}
If $R=E$, where $E$ is an infinite field, then the coefficient
polynomials are uniquely
determined by the equations (2.3) and the existence of such polynomials
for an $E{\rm GL}(n,E)$-module $U$ implies that $U$ is a polynomial module
of degree $r$.

References for basic facts about polynomial modules are only conveniently
available when $R$ is an infinite field. Thus the following result is often
useful.

\vspace{0.2in}
{\sc Lemma} 2.4. {\it Let\/ $U$ be an\/ $R{\rm GL}(n,R)$-module
with\/ $R$-basis\/ $\{u_1,\dots,u_d\}$. Let\/ $E$ be an infinite field
containing\/ $R$. Then\/ $U$ is a polynomial module of degree\/ $r$ if
and only if\/ $E \otimes_R U$ is a polynomial\/ $E{\rm GL}(n,E)$-module
which restricts to give\/ $U$, on taking\/ $U \subseteq
E \otimes_R U$ and ${\rm GL}(n,R) \subseteq {\rm GL}(n,E)$, and
such that the coefficient polynomials for\/ $E \otimes_R U$ with
respect to\/ $\{u_1,\dots,u_d\}$ belong to\/ $A_R(n,r)$.}

\vspace{0.2in}
{\it Proof}. (Sketch.) Suppose that $E \otimes_R U$ has the given
properties with coefficient polynomials $w_{st} \in A_R(n,r)$. By
(2.2), $U$ is invariant under the action of $S_R(n,r)$.
Conversely, if $U$ is an $S_R(n,r)$-module with coefficient
polynomials $w_{st}$ then (2.2) holds for $E \otimes_R U$ because
it holds for $U$. \hfill{$\square$}

\vspace{0.2in}
Let $V$ be a free $R$-module of rank $n$ with basis
$\{x_1,\dots,x_n\}$. The identity representation
${\rm GL}(n,R) \to {\rm GL}(n,R)$ gives $V$ the structure of the
``natural" $R{\rm GL}(n,R)$-module. It is easy to check (by Lemma 2.4,
for example) that $V$ is an $S_R(n,1)$-module. Indeed the
coefficient polynomials
are given by $w_{ij} = c_{ij}$ for all $i,j \in \{1,\dots,n\}$.
For every positive integer $r$, $T^r(V)$, $L^r(V)$ and $S^r(V)$ are
$R{\rm GL}(n,R)$-modules as we saw earlier in this section.
By Lemma 2.4, they are $S_R(n,r)$-modules. In \S3 we shall
consider $R{\rm GL}(n,R)$-modules of the form
$$S^{m(1)}(L^1(V)) \otimes_R S^{m(2)}(L^2(V)) \otimes_R \cdots
\otimes_R S^{m(r)}(L^r(V)),$$
where $m(1) + 2m(2) + \cdots + rm(r) = r$. Again they
are $S_R(n,r)$-modules.

From now on in this section we consider only the case $R = F$ where
$F$ is a field (finite or infinite). Thus $V$ is an $n$-dimensional
vector space over $F$.

For ${\bf i} \in I(n,r)$ where ${\bf i} = (i_1,\dots,i_r)$, write
$x_{\bf i} = x_{i_1}\cdots x_{i_r} \in T^r(V)$. Thus the elements
$x_{\bf i}$ form a basis of $T^r(V)$. The action of $S_F(n,r)$ on
$T^r(V)$ is given by [{\bf 16}, (2.6a)]. Translating to right modules we have
\begin{equation}
x_{\bf i} \xi = \sum_{\bf j} \xi(c_{{\bf i},\,{\bf j}}) x_{\bf j},
\end{equation}
for all $\xi \in S_F(n,r)$.

Let $\Lambda(n,r)$ be the set of all $n$-tuples
$\alpha = (\alpha_1,\dots,\alpha_n)$ of non-negative integers
satisfying $\alpha_1 + \cdots + \alpha_n = r$. For each such $\alpha$
take ${\bf i} \in I(n,r)$ so that
$c_{11}^{\alpha_1} c_{22}^{\alpha_2} \cdots c_{nn}^{\alpha_n} =
c_{{\bf i},{\bf i}}$
and write $\xi_\alpha = \xi_{{\bf i},{\bf i}}$. Thus (see [{\bf 16}, \S2.3
and \S3.2]) the elements $\xi_\alpha$ are pairwise orthogonal
idempotents of $S_F(n,r)$ which sum to the identity. It follows that if
$W$ is any $S_F(n,r)$-module then
\begin{equation}
W = \bigoplus_{\alpha \in \Lambda(n,r)} W\xi_\alpha,
\end{equation}
where the $W\xi_\alpha$ are subspaces, called weight spaces.

The multidegree of a monomial $x_{i_1}\cdots x_{i_r}$ of $T^r(V)$
is defined to be the element $\alpha = (\alpha_1,\dots,\alpha_n)$ of
$\Lambda(n,r)$ such that, for $j=1,\dots,n$, $\alpha_j$ is the
number of values of $t$ such that $i_t = j$. An element of $T^r(V)$
is said to be multihomogeneous if it is a linear combination of
monomials of the same multidegree.

By (2.4), the action of $\xi_\alpha$ on $T^r(V)$ is given by
\begin{equation}
(x_{i_1}\cdots x_{i_r})\xi_\alpha = \left\{
\begin{array}{ll} x_{i_1}\cdots x_{i_r} & \mbox{if $x_{i_1}\cdots x_{i_r}$ has
multidegree $\alpha$} \\
0 & \mbox{otherwise.}
\end{array} \right.
\end{equation}

Thus, if $W$ is an $S_F(n,r)$-submodule of $T^r(V)$, the weight space
$W\xi_\alpha$, as in (2.5),
consists of multihomogeneous elements of multidegree $\alpha$.

We shall need to consider $S_F(n,r)$-submodules of $T^r(V)$ as
$n$ varies. Thus, let $V^{(\omega)}$ be a vector space over $F$ with
countably infinite basis $\{x_1,x_2,\dots\}$. For each positive
integer $n$, let $V^{(n)}$ be the subspace with basis $\{x_1,\dots,x_n\}$ 
and identify $T(V^{(n)})$ with a subalgebra of
$T(V^{(\omega)})$ in the obvious way. Also, for each $n$, regard
$V^{(n)}$ as the natural $F{\rm GL}(n,F)$-module. Thus
$T^r(V^{(n)})$ is an $S_F(n,r)$-module.

Let $n_1$ and $n_2$ be positive integers with $n_1 \leq n_2$. Thus
$I(n_1,r) \subseteq I(n_2,r)$. For ${\bf i},{\bf j} \in I(n_1,r)$ we
may identify $c_{{\bf i},\,{\bf j}}$ in $A_F(n_1,r)$ with
$c_{{\bf i},\,{\bf j}}$ in $A_F(n_2,r)$ to make $A_F(n_1,r)$ a subspace
of $A_F(n_2,r)$. Similarly, by identifying elements $\xi_{{\bf i},\,{\bf j}}$,
we make $S_F(n_1,r)$ a subspace of $S_F(n_2,r)$. It is easy to verify
that $S_F(n_1,r)$ is a subalgebra of $S_F(n_2,r)$ (though the identity
elements are different if $n_1 < n_2$). We may also take
$\Lambda(n_1,r) \subseteq \Lambda(n_2,r)$ by identifying
$(\alpha_1,\dots,\alpha_{n_1})$ with $(\alpha_1,\dots,\alpha_{n_1},
0,\dots,0)$.

Let $W$ be an $S_F(n_2,r)$-submodule of $T^r(V^{(n_2)})$. Thus
$W = \bigoplus_{\alpha \in \Lambda(n_2,r)} W\xi_\alpha$. The
truncation function $d_{n_2,n_1}$ of [{\bf 16}, \S6.5] gives
$$d_{n_2,n_1}(W) = \bigoplus_{\alpha \in \Lambda(n_1,r)} W\xi_\alpha,$$
and $d_{n_2,n_1}(W)$ is an $S_F(n_1,r)$-submodule of
$T^r(V^{(n_1)})$. Another point of view is useful. Let
$\delta_{n_2,n_1}: T(V^{(n_2)}) \to T(V^{(n_1)})$ be the algebra
homomorphism defined by
$$\delta_{n_2,n_1}(x_i) = \left\{ \begin{array}{ll}
x_i & \mbox{for $i \in \{1,\dots,n_1\}$} \\
0 & \mbox{for $i \in \{n_1+1,\dots,n_2\}$.}
\end{array} \right. $$
Then it is easily verified, by (2.6), that
$$\delta_{n_2,n_1}(W) = d_{n_2,n_1}(W) = W \cap T(V^{(n_1)}).$$

Let $\{W^{(n)}: n \in {\mathbb N}\}$ be a family of subspaces of
$T(V^{(\omega)})$. (We write $\mathbb N$ for the set of all positive integers.)
We say that $\{W^{(n)}\}$ is a uniform family of modules of
degree $r$ if $W^{(n)}$ is an $S_F(n,r)$-submodule of $T^r(V^{(n)})$
for all $n$ and $\delta_{n_2,n_1}(W^{(n_2)}) = W^{(n_1)}$
for all $n_1,n_2$ with $n_2 \geq n_1$. For example, the family
$\{T^r(V^{(n)})\}$ is a uniform family of modules of degree $r$.

\vspace{0.2in}
{\sc Lemma} 2.5. {\it Let\/ $W$ be an\/ $S_F(r,r)$-submodule of\/
$T^r(V^{(r)})$. Then there exists a unique uniform family\/
$\{W^{(n)}\}$ of modules of degree\/ $r$ such that\/ $W^{(r)} = W$.}

\vspace{0.2in}
{\it Proof}. Suppose that $n \geq r$. We let the symmetric group
${\rm Sym}(n)$ act on
the right as a group of automorphisms of $T(V^{(n)})$ by extension
of its right
permutation action on $\{x_1,\dots,x_n\}$. For $\pi \in {\rm Sym}(n)$
and ${\bf i} \in I(n,r)$, where ${\bf i} = (i_1,\dots,i_r)$, we write
${\bf i}\pi = (i_1\pi,\dots,i_r\pi)$. Also, we identify $\pi$ with
the corresponding permutation matrix in ${\rm GL}(n,F)$ and write
$\tilde{\pi}$ for the element of $S_F(n,r)$ obtained from $\pi$
by the homomorphism (2.1). It is straightforward to verify that
in $S_F(n,r)$ we have
\begin{equation}
\tilde{\pi} \xi_{{\bf i},\,{\bf j}} = \xi_{{\bf i}\pi^{-1},\,{\bf j}}\;\;
{\rm and} \;\; \xi_{{\bf i},\,{\bf j}} \tilde{\pi} = \xi_{{\bf i},\,{\bf j}\pi},
\end{equation}
for all ${\bf i},{\bf j} \in I(n,r)$.

We shall show that $\sum_{\pi \in {\rm Sym}(n)} W\pi$ is an
$S_F(n,r)$-submodule of $T^r(V^{(n)})$. Let $\sigma \in {\rm Sym}(n)$ and
${\bf i},{\bf j} \in I(n,r)$. Choose $\pi \in {\rm Sym}(n)$ so that
${\bf j}\pi^{-1} \in I(r,r)$. By (2.7),
\begin{equation}
\tilde{\sigma} \xi_{{\bf i},\,{\bf j}} = \xi_{{\bf i}\sigma^{-1},\,{\bf j}} =
\xi_{{\bf i}\sigma^{-1},\,{\bf j}\pi^{-1}}\tilde{\pi}.
\end{equation}
For $w \in W$, $w\xi_{{\bf i}\sigma^{-1},\,{\bf j}\pi^{-1}} = 0$ if
${\bf i}\sigma^{-1} \notin I(r,r)$, by (2.4). However, if
${\bf i}\sigma^{-1} \in I(r,r)$, then
$\xi_{{\bf i}\sigma^{-1},\,{\bf j}\pi^{-1}} \in S_F(r,r)$, so
$w\xi_{{\bf i}\sigma^{-1},\,{\bf j}\pi^{-1}} \in W$. Thus, by (2.8),
$(W\sigma)\xi_{{\bf i},\,{\bf j}} \subseteq W\pi$. Hence
$\sum_{\pi \in {\rm Sym}(n)} W\pi$ is a submodule of $T^r(V^{(n)})$.

For $n \geq r$ we define $W^{(n)} = \sum_{\pi \in {\rm Sym}(n)} W\pi$
and for $n < r$ we define $W^{(n)} = \delta_{r,n}(W)$. Then it is
straightforward to verify that the family $\{W^{(n)}\}$ is unique
with the required properties. \hfill{$\square$}

\vspace{0.2in}
{\sc Lemma} 2.6.  {\it Let\/ $W_1$ and\/ $W_2$ be submodules of\/
$T^r(V^{(r)})$ which span their direct sum, and let\/ $\{W_1^{(n)}\}$,
$\{W_2^{(n)}\}$ and\/ $\{(W_1 \oplus W_2)^{(n)}\}$ be the uniform
families given by Lemma\/ {\rm 2.5}. Then, for all\/ $n$,
$$(W_1 \oplus W_2)^{(n)} = W_1^{(n)} \oplus W_2^{(n)}.$$}

\vspace{0.01in}
{\it Proof}. This is straightforward. \hfill{$\square$}

\vspace{0.2in}
{\sc Lemma} 2.7.  {\it Suppose that\/ $\{W_1^{(n)}\}$, \dots,
$\{W_m^{(n)}\}$ are uniform families of modules of
degrees\/ $r_1,\dots,r_m$,
respectively. Then\/
$\{\,[W_1^{(n)},W_2^{(n)},\dots,W_m^{(n)}]\,\}$
is a uniform family of modules of degree\/ $r_1 + \cdots + r_m$. }

\vspace{0.2in}
{\it Proof}. Let $r = r_1+\cdots +r_m$. By Lemma 2.4,
$[W_1^{(n)},\dots,W_m^{(n)}]$ is an $S_F(n,r)$-submodule
of $T^r(V^{(n)})$ for all $n$. Suppose that $n_1 \leq n_2$.
Then, since, $\delta_{n_2,n_1}$ is an algebra homomorphism,
\begin{eqnarray*}
\delta_{n_2,n_1}([W_1^{(n_2)},\dots,W_m^{(n_2)}]) & = &
[\delta_{n_2,n_1}(W_1^{(n_2)}),\dots,\delta_{n_2,n_1}(W_m^{(n_2)})] \\
& = & [W_1^{(n_1)},\dots,W_m^{(n_1)}].
\end{eqnarray*}
Thus the family is uniform. \hfill{$\square$}

\vspace{0.2in}
In particular, $\{L^r(V^{(n)}) : n \in {\mathbb N} \}$ is a uniform
family of modules of degree $r$.

We shall wish to consider arbitrary fields, possibly finite. However,
it is often convenient to pass to the infinite case, so we require
some results on field extensions.

\vspace{0.2in}
{\sc Lemma} 2.8. {\it Let\/ $F$ be a field and\/ $E$ an
extension field of\/ $F$. Let\/ $A$ be an\/ $F$-algebra and\/
$U$ a finite-dimensional\/ $A$-module. If\/ $E \otimes_F U$
is injective, as an\/ $E \otimes_F A$-module, then\/ $U$ is
injective. If\/ $E \otimes_F U$ is projective then\/ $U$ is
projective.}

\vspace{0.2in}
{\it Proof}. Suppose that $E \otimes U$ is injective as
an $E \otimes A$-module. Then, by [{\bf 11}, Proposition II.6.2a],
$E \otimes U$ is injective as an $A$-module. But $U$ is a direct
summand of $E \otimes U$ as an $A$-module. Hence $U$ is injective.
The result for projective modules is proved in the same way, using
[{\bf 11}, Proposition II.6.2]. \hfill{$\square$}

\vspace{0.2in}
As above, we consider $V^{(n)}$ over an arbitrary field $F$.

\vspace{0.2in}
{\sc Lemma} 2.9. {\it Suppose
that\/ $r \leq n$ and let\/ $U$ be an\/ $S_F(n,r)$-module which
is isomorphic to a direct summand of\/ $T^r(V^{(n)})$. Then\/ $U$ is
injective and projective as an\/ $S_F(n,r)$-module.}

\vspace{0.2in}
{\it Proof}. Let $E$ be an infinite extension field of $F$. Thus
$E \otimes U$ is isomorphic to a direct summand of
$T^r(E \otimes V^{(n)})$.
Since $r \leq n$, $T^r(E \otimes V^{(n)})$ is injective as an
$S_E(n,r)$-module, by [{\bf 12}, (3.4) Lemma].
However, by [{\bf 16}, Example 1 after (2.7e)], $T^r(E \otimes V^{(n)})$
is isomorphic to its contravariant dual. Thus $T^r(E \otimes V^{(n)})$
is projective. It follows that $E \otimes U$ is injective and
projective. Therefore, by Lemma 2.8, $U$ is injective and
projective. \hfill{$\square$}

\vspace{0.2in}
For positive integers $n$ and $r$, let ${\cal T}_r^{(n)}$ denote the
class of all finite-dimensional $S_F(n,r)$-modules that are direct sums
of modules each of which is isomorphic to a direct summand of
$T^r(V^{(n)})$.

\vspace{0.2in}
{\sc Lemma} 2.10. {\it Let\/ $\{W^{(n)}\}$ be a uniform family of
modules of degree\/ $r$ such that\/ $W^{(r)} \in {\cal T}_r^{(r)}$.
Then, for all\/ $n$, $W^{(n)}$ is a direct summand
of\/ $T^r(V^{(n)})$.}

\vspace{0.2in}
{\it Proof}. By Lemma 2.9, $W^{(r)}$ is injective. Therefore, since
$W^{(r)}$ is a submodule of $T^r(V^{(r)})$, it is a
direct summand of $T^r(V^{(r)})$. It follows, by Lemma 2.6, that
$W^{(n)}$ is a direct summand of $T^r(V^{(n)})$ for all $n$.
\hfill{$\square$}

\vspace{0.4in}
\begin{center}
3. {\it Filtrations of tensor powers}
\end{center}
\renewcommand{\theequation}{3.\arabic{equation}}
\setcounter{equation}{0}

As in \S2, we first work over a ring $R$ which is either the
ring of integers or an arbitrary field.  Let $V$ be a free $R$-module
and let $r$ be a positive integer.

By a composition of $r$ we mean a sequence $\lambda = (\lambda_1, \dots,
\lambda_l)$ of positive integers $\lambda_1,\dots,\lambda_l$
satisfying $\lambda_1 + \cdots + \lambda_l = r$. If we also have
$\lambda_1 \geq \cdots \geq \lambda_l$ then, as usual, $\lambda$ is called a
partition of $r$.
If $\mu = (\mu_1,\dots,\mu_l)$ is a composition of $r$
and we put $\mu_1,\dots,\mu_l$ into non-increasing order we obtain
a partition of $r$, which we call the partition associated to $\mu$.
We write ${\rm Part}(r)$ for the set of all partitions of
$r$.

If $\lambda, \theta \in {\rm Part}(r)$, where
$\lambda = (\lambda_1,\dots,\lambda_k)$ and $\theta =
(\theta_1,\dots,\theta_l)$, we say that $\lambda$ is a
refinement of $\theta$ if there is a function
$f: \{1,\dots,k\} \to \{1,\dots,l\}$ such that
$\theta_j = \sum_{i:f(i) = j}\lambda_i$ for $j=1,\dots,l$.
For example, $(3,2,2,2,1)$ is a refinement of $(4,3,3)$.
The relation of refinement is a partial order on ${\rm Part}(r)$,
and this partial order plays a key role in the study of
filtrations of tensor powers.
However, it is more convenient here to work with a total order
$\leq^*$ which extends refinement.
We use the ``lexicographic" order in which
$$(\lambda_1,\lambda_2,\dots,\lambda_k) <^* (\theta_1,
\theta_2, \dots,\theta_l)$$
if there exists $i$ such that $\lambda_i < \theta_i$ but
$\lambda_j = \theta_j$ for all $j < i$. Note that
if $\lambda$ is a refinement of $\theta$ then $\lambda \leq^* \theta$.
Note also that the smallest partition is $(1,1,\dots,1)$, written as 
$(1^r)$, and the largest partition is $(r)$.

Our main results are stated in terms of the lexicographic order
$\leq^*$,
but similar results apply for any total order which extends refinement, 
and more complicated versions of the results could be stated in terms of
refinement itself.

Let ${\cal X}_r$ be the set of all elements of $T^r(V)$ of the
form $b_1b_2\cdots b_l$ ($l$ arbitrary) where each $b_i$ belongs
to $L^{\mu_i}(V)$ for some positive integer $\mu_i$. Thus $\mu_i = \deg b_i$
for all $i$ and $(\mu_1,\dots,\mu_l)$ is a composition of $r$. For
each partition $\lambda$ of $r$ let ${\cal X}_\lambda$ denote the
set of all such elements $b_1b_2\cdots b_l$ where
$(\deg b_1,\dots,\deg b_l)$ has $\lambda$ as its associated
partition.

For each $\lambda$,  let $W_\lambda$ be the
$R$-submodule of $T^r(V)$
spanned by all ${\cal X}_{\theta}$ with $\lambda 
\leq^* \theta$. For each $\lambda$
such that $\lambda \neq (r)$ let $\lambda+$ be the partition of $r$
which is next bigger than $\lambda$. Thus we have the filtration
$$T^r(V) = W_{(1^r)} \geq \cdots \geq
W_\lambda \geq W_{\lambda+} \geq \cdots \geq W_{(r)} \geq 0.$$
We call this the PBW-filtration of $T^r(V)$. (As we shall
see it is closely connected with the Poincar\'e--Birkhoff--Witt
Theorem.) It is convenient to define $W_{(r)+} = 0$. Thus, for all
$\lambda$, ${\cal X}_\lambda$ spans $W_\lambda$ modulo $W_{\lambda+}$.

The following result is well known and easy to prove using the identity
$b_ib_j = b_jb_i + [b_i,b_j]$.

\vspace{0.2in}
{\sc Lemma} 3.1. {\it Let\/ $\mu = (\mu_1,\dots,\mu_l)$ be a
composition of\/ $r$ with associated partition\/ $\lambda$.
For\/ $i=1,\dots,l$, let\/ $b_i \in L^{\mu_i}(V)$, and let\/
$\pi \in {\rm Sym}(l)$: thus\/ $b_1b_2\cdots b_l$ and\/
$b_{1\pi}b_{2\pi}\cdots b_{l\pi}$ belong to\/ ${\cal X}_\lambda$.
We have
$$b_1b_2\cdots b_l + W_{\lambda+} = b_{1\pi}b_{2\pi}\cdots b_{l\pi}
+ W_{\lambda+}.$$}

\vspace{0.01in}
For each positive
integer $i$ choose an $R$-basis ${\cal A}_i$ of $L^i(V)$.
The union of these is
an $R$-basis $\cal A$ of $L(V)$. Take any total ordering
$\preccurlyeq$ of $\cal A$. Then,
by the Poincar\'e--Birkhoff--Witt Theorem (see [{\bf 18}] for the case of a
field and [{\bf 2}] for general $R$), $T(V)$ has an
$R$-basis $\cal F$ consisting of all elements of the form
\begin{equation}
a_1a_2\cdots a_l \quad (l \geq 0, \; a_1,\dots,a_l \in {\cal A}, \;
a_1 \preccurlyeq a_2 \preccurlyeq \cdots \preccurlyeq a_l).
\end{equation}
We call $\cal F$ a PBW-basis of $T(V)$. Note that the elements (3.1)
of $\cal F$ are distinct as written.
We write ${\cal F}_r$ for the set of all elements (3.1) of degree $r$.
Also, for each $\lambda \in {\rm Part}(r)$,
we write ${\cal F}_\lambda$ for the set of elements (3.1) in ${\cal F}_r$
such that the composition $(\deg a_1, \dots, \deg a_l)$
has $\lambda$ as associated partition. Clearly, ${\cal F}_r$ is an
$R$-basis of $T^r(V)$, called a PBW-basis of $T^r(V)$. Also, it is
easy to verify, by Lemma 3.1 and reverse induction on $\lambda$
(starting with $\lambda = (r)$), that the
elements of ${\cal F}_\lambda$ taken modulo $W_{\lambda+}$ form
an $R$-basis of $W_\lambda/W_{\lambda+}$. (It follows that $W_\lambda >
W_{\lambda +}$ provided that $\dim V > 1$.)

Suppose that $V$ is an $RG$-module for some group $G$. Let $\lambda \in
{\rm Part}(r)$ and write $\lambda$ in the form
$$\lambda = (r^{m(r)}, (r-1)^{m(r-1)}, \dots, 1^{m(1)}),$$
where $m(r),\dots, m(1)$ are non-negative integers. Then we define an
$RG$-module $L^\lambda(V)$ by
$$L^\lambda(V) =
S^{m(1)}(L^1(V)) \otimes_R S^{m(2)}(L^2(V)) \otimes_R \cdots
\otimes_R S^{m(r)}(L^r(V)).$$
This is called the higher Lie module corresponding to $\lambda$.

The following result is essentially well known: see, for example,
[{\bf 17}, Lemma 3.3]. However, since we have no fully appropriate
reference we give a sketch proof.

\vspace{0.2in}
{\sc Lemma} 3.2. {\it Suppose that the free\/ $R$-module\/ $V$ is
an\/ $RG$-module for some group\/ $G$. Then,
for each\/ $\lambda \in {\rm Part}(r)$,
$W_\lambda$ is an\/ $RG$-submodule of\/ $T^r(V)$ and\/
$W_\lambda/W_{\lambda+} \cong L^\lambda(V)$.}

\vspace{0.2in}
{\it Proof}. (Sketch.) It is easily verified by reverse induction on
$\lambda$ that $W_\lambda$ is an $RG$-submodule of $T^r(V)$.
In the notation used above, choose the ordering of
$\cal A$ in any way so that all elements of ${\cal A}_i$
are smaller than all elements of ${\cal A}_j$ whenever $i < j$.
The basis $\{w + W_{\lambda+} : w \in {\cal F}_\lambda \}$ of
$W_\lambda/W_{\lambda+}$ consists of all elements of the form
$$a_1^{(1)}\cdots a_1^{(m(1))}a_2^{(1)}\cdots a_2^{(m(2))}
\cdots a_r^{(1)}\cdots a_r^{(m(r))} + W_{\lambda+}$$
where the factors are in non-decreasing order and
$a_i^{(1)},\dots,a_i^{(m(i))} \in {\cal A}_i$ for $i=1,\dots,r$.
There is a corresponding basis of $L^\lambda(V)$
consisting of the elements
$$a_1^{(1)} \cdots a_1^{(m(1))} \otimes
a_2^{(1)} \cdots a_2^{(m(2))} \otimes \cdots \otimes
a_r^{(1)} \cdots a_r^{(m(r))}.$$
It is easy to verify, by Lemma 3.1, that the
$R$-module isomorphism $\alpha$ that takes basis element to corresponding
basis element is an $RG$-module isomorphism. \hfill{$\square$}

\vspace{0.2in}
In particular, $W_{(1^r)}/W_{(1^r)+} \cong S^r(V)$ and
$W_{(r)} \cong L^r(V)$. 

If $V$ has finite rank $n$ then the isomorphism of
Lemma 3.2 is an $S_R(n,r)$-module isomorphism. This is clear if
$R$ is an infinite field (because then the homomorphism (2.1) is
surjective). In the general case, let $E$ be an
infinite field containing $R$. By Lemma 3.2, $W_\lambda$ is
an $R{\rm GL}(n,R)$-submodule of $T^r(V)$. We regard $T^r(V)$ as a
subset of $T^r(E \otimes V)$ and identify this with $E \otimes T^r(V)$.
It is easily verified that $E \otimes W_\lambda$ is the subspace
corresponding to $\lambda$ in the PBW-filtration of $T^r(E \otimes V)$.
Thus we may verify by Lemma 2.4 that $W_\lambda$ is an
$S_R(n,r)$-submodule of $T^r(V)$. Let $\alpha$ be the
$R$-module isomorphism from $W_\lambda/W_{\lambda+}$ to
$L^\lambda(V)$ described in the proof of Lemma 3.2.
Upon extension to $E$, $\alpha$
yields an $S_E(n,r)$-module isomorphism, by the case first considered.
Therefore $\alpha$ is an $S_R(n,r)$-module isomorphism.

\vspace{0.2in}
We shall make use of Solomon's descent algebra. For definitions and further
details see [{\bf 14}, \S2] and the works cited there.
As well as results from [{\bf 14}] we shall use a theorem of Garsia and
Reutenauer [{\bf 15}]. Our notation is
similar to that in [{\bf 14}, \S2]. Let $D_r$ be the descent algebra of
${\rm Sym}(r)$ over the integers. Thus $D_r$ is a subalgebra of
${\mathbb Z}{\rm Sym}(r)$ and
$D_r$ has a $\mathbb Z$-basis consisting of certain elements
$X^\nu$ of ${\mathbb Z}{\rm Sym}(r)$ indexed by the compositions $\nu$
of $r$. (In the notation of [{\bf 15}], $X^\nu$ would be written as $B_\nu$.)
The elements $X^\nu$ can be interpreted as elements of
$R\,{\rm Sym}(r)$ and these elements
form an $R$-basis for an $R$-subalgebra $D_{r,R}$ of $R\,{\rm Sym}(r)$.

For each composition $\nu$ of $r$, where $\nu = (\nu_1,\dots,\nu_k)$,
we write $\phi^\nu$ for the corresponding Young character of
${\rm Sym}(r)$. This is the character (in characteristic $0$) induced
from the trivial character of the group ${\rm Sym}(\nu_1) \times \cdots
\times {\rm Sym}(\nu_k)$ regarded as a subgroup of ${\rm Sym}(r)$ in the
obvious way.  The $\mathbb Z$-span of the Young characters is a subring
of the ring of all class functions ${\rm Sym}(r) \to {\mathbb Z}$.
For each $\nu$, define $\phi^{\nu,R} : {\rm Sym}(r) \to R$ by
$\phi^{\nu,R}(\sigma) = \phi^\nu(\sigma)1_R$ for all
$\sigma \in {\rm Sym}(r)$. Thus the $R$-span of the
$\phi^{\nu,R}$ is a subalgebra $C_{r,R}$ of the $R$-algebra of all
class functions ${\rm Sym}(r) \to R$. 
For $\xi \in C_{r,R}$ and
$\lambda \in {\rm Part}(r)$ we write $\xi(\lambda)$ for the value of
$\xi$ at the conjugacy class of ${\rm Sym}(r)$ indexed by $\lambda$.
By the same proof as for [{\bf 14}, Theorem 3], there is a surjective
homomorphism of algebras $c_{r,R} : D_{r,R} \to C_{r,R}$ satisfying
$c_{r,R}(X^\nu) = \phi^{\nu,R}$ for all compositions $\nu$ of $r$.

The symmetric group ${\rm Sym}(r)$ has a left action on $T^r(V)$ by
``place permutations". This action is given by
$\sigma(v_1v_2\cdots v_r) = v_{1\sigma}v_{2\sigma} \cdots
v_{r\sigma}$
for all $\sigma \in {\rm Sym}(r)$, $v_1,\dots,v_r \in V$. Thus
$T^r(V)$ is an $R\,{\rm Sym}(r)$-module and hence a $D_{r,R}$-module.

\vspace{0.2in}
{\sc Proposition} 3.3. {\it For each partition\/ $\lambda$ of\/ $r$,
$W_\lambda$ is a\/ $D_{r,R}$-submodule of\/ $T^r(V)$. For all\/
$X \in D_{r,R}$, the action of\/ $X$ on\/ $W_\lambda/W_{\lambda+}$
is given by multiplication by the element\/ $(c_{r,R}(X))(\lambda)$
of\/ $R$.}

\vspace{0.2in}
{\it Proof}.
Let $\mu$ and $\nu$ be compositions of $r$, where
$\mu = (\mu_1,\dots,\mu_l)$ and $\nu = (\nu_1,\dots,\nu_k)$.
Let $Q(\nu,\mu)$ be the set of all $k$-tuples
$(I_1,\dots,I_k)$ where $I_1,\dots,I_k$ are subsets of $\{1,\dots,l\}$
such that $\{1,\dots,l\}$ is the disjoint union of $I_1,\dots,I_k$ and
$\nu_j = \sum_{i \in I_j}\mu_i$ for $j=1,\dots,k$.
It is easy to verify that $|Q(\nu,\mu)|$ is the number of $k \times l$
matrices of non-negative integers such that (the $l$-tuple formed by)
the sum of the rows is $\mu$, the sum of the columns is $\nu$ and
each column contains a unique non-zero entry. Therefore, by [{\bf 1}, (12)],
\begin{equation}
|Q(\nu,\mu)| = \phi^\nu(\lambda),
\end{equation}
where $\lambda$ is the partition associated to $\mu$. 

For each
$q = (I_1,\dots,I_k) \in Q(\nu,\mu)$ let $\pi_q$ be the permutation
of $\{1,\dots,l\}$ such that $(1\pi_q,2\pi_q,\dots,l\pi_q)$ consists
of the elements of $I_1$ in increasing order followed by the elements
of $I_2$ in increasing order, and so on. For $i=1,\dots,l$, let
$b_i \in L^{\mu_i}(V)$. Then, by [{\bf 15}, Theorem 2.1],
\begin{equation}
X^\nu(b_1b_2\cdots b_l) = \sum_{q \in Q(\nu,\mu)}
b_{1\pi_q}b_{2\pi_q}\cdots b_{l\pi_q}.
\end{equation}
Note that this result, as stated in [{\bf 15}], applies only to the case
where $R = {\mathbb Q}$ and $V$ has finite rank. However, since (3.3)
involves only finitely many elements of $V$, the same result holds for
$V$ of infinite rank. From the result over ${\mathbb Q}$ we can deduce
the result over ${\mathbb Z}$ because, if $V$ is a
free $\mathbb Z$-module, there are natural embeddings of $T(V)$
into $T({\mathbb Q} \otimes V)$ and of $D_r$ into $D_{r,{\mathbb Q}}$.
In the case where $R$ is an arbitrary field, let $V$ have basis
$\{x_1,x_2,\dots\}$. Both sides of (3.3) are linear in each term $b_i$.
Thus it suffices to prove (3.3) in the case where each $b_i$ is a Lie
monomial in $x_1,x_2,\dots\;$. The calculation of $X^\nu(b_1b_2\cdots b_l)$
is then the same as it is over ${\mathbb Z}$
followed by reduction modulo the
characteristic of $R$. Thus (3.3) holds in all cases.

Let $\lambda \in {\rm Part}(r)$. Then,
by (3.3), we have $X^\nu u \in W_\lambda$ for all
$u \in {\cal X}_\lambda$ and all compositions $\nu$ of $r$. Thus
$Xu \in W_\lambda$ for all $X \in D_{r,R}$.
It follows that each $W_\lambda$ is
a $D_{r,R}$-submodule of $T^r(V)$.  Hence $W_\lambda/W_{\lambda+}$ is a
$D_{r,R}$-module. Let $b_1b_2 \cdots b_l \in {\cal X}_\lambda$, with
composition $\mu = (\deg b_1,\dots,\deg b_l)$, and
let $\nu$ be any composition of $r$.
By (3.3) and Lemma 3.1,
$$X^\nu(b_1b_2\cdots b_l) + W_{\lambda+} = |Q(\nu,\mu)| b_1b_2\cdots b_l
+ W_{\lambda+}.$$
Thus, by (3.2), the action of $X^\nu$ on $W_\lambda/W_{\lambda+}$ is
multiplication by the scalar $\phi^\nu(\lambda)1_R$. By the definition of
$c_{r,R}$, this scalar is $(c_{r,R}(X^\nu))(\lambda)$. This holds for
all $\nu$, giving the required result. \hfill{$\square$}

\vspace{0.2in}
From now on in this section we assume that $R=F$, where $F$ is a field
of prime characteristic $p$.
If $\lambda,\theta \in {\rm Part}(r)$ and $\sigma_\lambda$ and
$\sigma_\theta$
are elements of ${\rm Sym}(r)$ in the conjugacy classes corresponding to
$\lambda$ and $\theta$, respectively, then,
as in [{\bf 14}, \S2], we
say that $\lambda$ and $\theta$ are $p$-equivalent
if $\sigma_\lambda^{p^s}$ is
conjugate to $\sigma_\theta^{p^s}$ for some positive integer $s$. This
gives an equivalence relation on ${\rm Part}(r)$. By [{\bf 14},
Proposition 5], for all $\xi \in C_{r,F}$ we have
$\xi(\lambda) = \xi(\theta)$ whenever $\lambda$
and $\theta$ are $p$-equivalent.

Let $\cal J$ be the set of all $p$-equivalence classes of partitions
of $r$. As shown in [{\bf 14}, Corollary 6],
the identity element $1$ of $D_{r,F}$ can be written as
$1 = \sum_{J \in {\cal J}} e_J$, where the $e_J$ are
pairwise-orthogonal primitive idempotents in $D_{r,F}$ 
such that $c_{r,F}(e_J)$ is given by
$$(c_{r,F}(e_J))(\lambda) = \left\{ \begin{array}{ll}
1 & \mbox{if $\lambda \in J$} \\
0 & \mbox{if $\lambda \notin J$.} \end{array} \right. $$
Therefore, by Proposition 3.3, $e_J$ acts on
$W_\lambda/W_{\lambda+}$ as the identity if
$\lambda \in J$ and as zero otherwise. Thus
$$e_JW_\lambda + W_{\lambda+} = \left\{ \begin{array}{ll}
W_\lambda & \mbox{if $\lambda \in J$} \\
W_{\lambda+} & \mbox{ if $\lambda \notin J$.} \end{array} \right.$$
Also,
\begin{equation}
e_Jw + W_{\lambda+} = w + W_{\lambda+} \quad \mbox{for all
$w \in W_\lambda$, where $\lambda \in J$.}
\end{equation}
Clearly we obtain a chain of subspaces
\begin{equation}
e_JT^r(V) = e_JW_{(1^r)} \geq \cdots \geq e_JW_\lambda \geq
e_JW_{\lambda+} \geq \cdots \geq e_JW_{(r)} \geq 0.
\end{equation}
Note that, since $e_J$ is an
idempotent, we have $e_JW_\lambda \cap W_{\lambda+} = e_JW_{\lambda+}$.
Thus $e_JW_\lambda/e_JW_{\lambda+} \cong (e_JW_\lambda + W_{\lambda+})
/W_{\lambda+}$ and we have
\begin{equation}
\mbox{$e_JW_\lambda/e_JW_{\lambda+} \cong
W_\lambda/W_{\lambda+}\;$
if $\lambda \in J$,}
\end{equation}
while
\begin{equation}
\mbox{$e_JW_\lambda/e_JW_{\lambda+} = 0\;$ if $\lambda \notin J$.}
\end{equation}
Therefore, if $J = \{\lambda(1),\dots,\lambda(s)\}$ with
$\lambda(1) <^* \cdots <^* \lambda(s)$, we may write the chain (3.5) in
the form
$$e_JT^r(V) = U_1 \geq U_2 \geq \cdots \geq U_s \geq U_{s+1} = 0,$$
where $U_j/U_{j+1} \cong W_{\lambda(j)}/W_{\lambda(j)+}$ for
$j=1,\dots,s$.

Suppose that $G$ is a group and $V$ is an $FG$-module.
It is well known and easy to verify that
the actions of ${\rm Sym}(r)$ and $G$ on $T^r(V)$ commute.
It follows that $e_JT^r(V)$ is an $FG$-submodule of $T^r(V)$ and the
function $\phi_J : T^r(V) \to e_JT^r(V)$ given by $w \mapsto e_Jw$ for all
$w \in T^r(V)$ is a homomorphism of $FG$-modules. By the properties
of the idempotents $e_J$ we see that $T^r(V)$ is the direct sum of the
submodules $e_JT^r(V)$. The chain (3.5) is a filtration of
$e_JT^r(V)$ by $FG$-submodules. Also, for $\lambda \in J$, the
isomorphism (3.6) is evidently an isomorphism of $FG$-modules.
In view of Lemma 3.2, we have proved the following result.

\vspace{0.2in}
{\sc Theorem} 3.4 ({\sc Filtration Theorem}).
{\it Let\/ $F$ be a field of prime characteristic\/ $p$.
Let\/ $G$ be a group and\/ $V$ an\/ $FG$-module.
Let\/ $e_J$, $J \in {\cal J}$, be primitive idempotents
of the descent algebra as described above.
Then the\/ $FG$-module\/
$T^r(V)$ is the direct sum of its submodules\/ $e_JT^r(V)$.
For a given\/ $J$, let\/ $\lambda(1)$, \dots, $\lambda(s)$ be
the partitions belonging to\/ $J$ in increasing lexicographic
order. Then\/ $e_JT^r(V)$ has a filtration of\/ $FG$-submodules,
$$e_JT^r(V) = U_1 \geq U_2 \geq \cdots \geq U_s \geq U_{s+1} = 0,$$
where, for\/ $j=1,\dots,s$, we have\/ $U_j/U_{j+1}
\cong L^{\lambda(j)}(V)$. }

\vspace{0.2in}
It is easily seen that if $V$ has finite dimension $n$ then each of the
isomorphisms of Theorem 3.4 is an isomorphism of $S_F(n,r)$-modules.

Some supplementary details will be important later, 
so we note them in the proposition below. For each $\lambda
\in {\rm Part}(r)$, let ${\cal Y}_\lambda$ be any subset of $W_\lambda$
such that the family $\{y + W_{\lambda+} : y \in {\cal Y}_\lambda\}$
consists of distinct elements forming a
basis of $W_\lambda/W_{\lambda+}$. Thus $\bigcup {\cal Y}_\lambda$ is
a basis of $T^r(V)$. For each $\lambda$ let
$Y_\lambda$ be the subspace of $T^r(V)$ with basis ${\cal Y}_\lambda$.
Let $J$ be a $p$-equivalence class of partitions of $r$.
Let ${\cal Y}_J = \bigcup_{\lambda \in J} {\cal Y}_\lambda$ and let
$Y_J$ be the subspace of $T^r(V)$ with basis ${\cal Y}_J$.
Thus $Y_J = \bigoplus_{\lambda \in J}Y_\lambda$.

\vspace{0.2in}
{\sc Proposition} 3.5. {\it In the notation above, the function\/
$\phi_J: T^r(V) \to
e_JT^r(V)$ given by\/ $\phi_J(w) = e_Jw$ for all\/ $w \in T^r(V)$
restricts to a vector space isomorphism from\/ $Y_J$ to\/
$e_JT^r(V)$.  Thus\/ $e_J{\cal Y}_J$ is a basis of\/ $e_JT^r(V)$.}

\vspace{0.2in}
{\it Proof}.  By (3.7), we have
$$\mbox{$\phi_J(W_\lambda) =
\phi_J(W_{\lambda+})\;$ if $\lambda \notin J$.}$$
Also, since $W_\lambda = Y_\lambda + W_{\lambda+}$, we have
$$\mbox{$\phi_J(W_\lambda) \subseteq \phi_J(Y_J) +
\phi_J(W_{\lambda+})\;$ if $\lambda \in J$.}$$
It follows, by reverse induction on $\lambda$ (starting with
$\lambda = (r)$), that $\phi_J(W_\lambda) \subseteq \phi_J(Y_J)$
for all $\lambda$. Therefore $\phi_J(T^r(V)) \subseteq \phi_J(Y_J)$
and so $\phi_J(Y_J) = \phi_J(T^r(V)) = e_JT^r(V)$. Thus the restriction
$\phi_J:Y_J \to e_JT^r(V)$ is surjective. Suppose that $y$ is a non-zero
element of $Y_J$. Then, for some $\lambda \in J$, we have
$y \in W_\lambda \smallsetminus W_{\lambda+}$. Since $e_J$ acts
as the identity on $W_\lambda/W_{\lambda+}$, we have
$\phi_J(y) \in W_\lambda \smallsetminus W_{\lambda+}$. Thus
$\phi_J(y) \neq 0$ and so $\phi_J:Y_J \to e_JT^r(V)$ is injective.
\hfill{$\square$}

\vspace{0.4in}
\begin{center}
4. {\it The Decomposition Theorem}
\end{center}
\renewcommand{\theequation}{4.\arabic{equation}}
\setcounter{equation}{0}

We begin with a major step towards the proof.

\vspace{0.2in}
{\sc Proposition} 4.1. {\it
Let\/ $V$ be a finite-dimensional vector space over a field\/ $F$
of prime characteristic\/ $p$. Regard\/ $V$ as the natural\/
$F{\rm GL}(n,F)$-module, where\/ $n = \dim V$. Let\/ $q$ be a
positive integer such that\/ $q \leq n$ and write\/
$q = p^mq'$ where\/ $q'$ is not divisible by\/ $p$. For each
proper divisor\/ $d$ of\/ $q$, let\/ $B_d$ be an\/ $S_F(n,d)$-submodule
of\/ $L^d(V)$ such that\/ $B_d$ is a direct summand of\/ $T^d(V)$.
Suppose that 
$$L^{p^iq'}(V) = \bigoplus_{d \mid p^iq'} L^{p^iq'/d}(B_d),$$
for\/ $i=0,\dots,m-1$, 
and suppose that the modules\/ $L^{q/d}(B_d)$, where\/ $d$ runs
over the proper divisors of\/ $q$, span their direct sum\/ $C$
in\/ $L^q(V)$. Then\/ $L^q(V)/C$ is injective and projective
as an\/ $S_F(n,q)$-module.}

\vspace{0.2in}
Note that the notation $L^r(B_d)$ is justified by the remarks
early in \S2.

\vspace{0.2in}
{\it Proof.}
We first reduce to the case where $F$ is
infinite. Let $E$ be an infinite extension field of $F$ and
identify $E \otimes T(V)$ with $T(E \otimes V)$. Thus
$E \otimes L^r(V)$ is identified with $L^r(E \otimes V)$ for all $r$.
The modules $E \otimes B_d$ within $T(E \otimes V)$ satisfy
the same hypotheses as the modules $B_d$ within $T(V)$.
Suppose that the proposition holds over $E$. Then
$E \otimes L^q(V)/E \otimes C$ is injective and projective as
an $S_E(n,q)$-module. However, this module is isomorphic
to $E \otimes (L^q(V)/C)$ and we may identify
$S_E(n,q)$ with $E \otimes S_F(n,q)$. Therefore, by Lemma 2.8,
$L^q(V)/C$ is injective and projective as an $S_F(n,q)$-module.

Therefore we may assume that $F$ is infinite. The advantage is that
$F{\rm GL}(n,F)$-submodules of $T^r(V)$ are the same as
$S_F(n,r)$-submodules.

If $q = 1$ then $C = 0$ and $L^q(V) = V$: thus the required
conclusion follows from Lemma 2.9. Hence we may assume that $q > 1$.

The remaining proof will be done in three steps. Taking careful
account of the hypotheses, we first choose an ordered basis
$\cal A$ of $L(V)$ which contains a basis of $C$. This yields a
PBW-basis ${\cal F}_q$ of $T^q(V)$. We then use idempotents from
various descent algebras to modify ${\cal F}_q$ to another
basis ${\cal F}_q^*$ of $T^q(V)$. Certain elements of this
modified basis will be shown to span a direct summand $U$ of
$T^q(V)$ such that $U \cap L^q(V) = C$ and such that $U + L^q(V)$
is also a direct summand of $T^q(V)$. The claim about $L^q(V)/C$
then follows easily.

Let $\bf D$ denote the set of all proper divisors of $q$. For each
$d \in {\bf D}$ we write $d = p^{j(d)}d'$ where $d'$ is not
divisible by $p$. For $d \in {\bf D}$ and $i \in \{0,\dots,m\}$
such that $d \mid p^iq'$ (equivalently, $i \in \{j(d),\dots,m\}$),
choose a basis ${\cal B}(p^iq'/d,d)$ of $L^{p^iq'/d}(B_d)$.

For $i \in \{0,\dots,m-1\}$ write
$${\cal A}_{p^iq'} = \bigcup_{d \mid p^iq'} {\cal B}(p^iq'/d,d).$$
Thus, by the hypotheses of the proposition, ${\cal A}_{p^iq'}$
is a basis of $L^{p^iq'}(V)$. Similarly, $\bigcup_{d \in {\bf D}}
{\cal B}(q/d,d)$ is a basis of $C$. We extend this to a basis
${\cal A}_q$ of $L^q(V)$: thus
$${\cal A}_q = \bigcup_{d \in {\bf D}} {\cal B}(q/d,d)
\cup {\cal A}_q',$$
for some complementary set ${\cal A}_q'$. For values of $r$ not
already dealt with, let ${\cal A}_r$ be an arbitrary basis of
$L^r(V)$, and write ${\cal A} = \bigcup_{r \in {\mathbb N}}
{\cal A}_r$. Thus $\cal A$ is a basis of $L(V)$.

For each $d \in {\bf D}$ write
\begin{equation}
{\cal B}(d) = \bigcup_{i \in \{j(d),\dots,m\}}
{\cal B}(p^iq'/d,d).
\end{equation}
Thus
\begin{equation}
\bigcup_{d \in {\bf D}} {\cal B}(d) =
\biggl(\bigcup_{i \in \{0,\dots,m\}} {\cal A}_{p^iq'} \biggr)
\smallsetminus {\cal A}_q'.
\end{equation}

We order ${\cal A}$ in any way subject to all elements of
${\cal B}(d_1)$ being smaller than all elements of
${\cal B}(d_2)$ whenever $d_1,d_2 \in {\bf D}$ with $d_1 < d_2$.
As in \S3, we obtain a PBW-basis ${\cal F}_q$ of $T^q(V)$
consisting of all elements of the form $a_1a_2\cdots a_l$, with
$l \geq 0$, $a_1,\dots,a_l \in {\cal A}$, $a_1 \preccurlyeq \cdots
\preccurlyeq a_l$
and $\deg a_1\cdots a_l = q$.
For each $\lambda \in {\rm Part}(q)$, let ${\cal F}_\lambda$
and $W_\lambda$ be defined as in \S3.
Thus we have the PBW-filtration
$$T^q(V) = W_{(1^q)} \geq \cdots
\geq W_\lambda \geq W_{\lambda+} \geq \cdots \geq W_{(q)} \geq 0,$$
and the elements of
${\cal F}_\lambda$ taken modulo $W_{\lambda+}$ form a basis of
$W_\lambda/W_{\lambda+}$.

If we replace each element
$a_1a_2\cdots a_l$ of ${\cal F}_\lambda$ by any element of $W_\lambda$
which is equal to $a_1a_2\cdots a_l$ modulo
$W_{\lambda+}$ we get a new set ${\cal F}^*_\lambda$ and if we do this
for each $\lambda$ we get a new basis ${\cal F}_q^* = \bigcup
{\cal F}^*_\lambda$ of $T^q(V)$.

For each positive integer $r$ and each $\lambda \in {\rm Part}(r)$
let ${\rm cl}(\lambda)$ denote the
$p$-equivalence class of $\lambda$ in ${\rm Part}(r)$.
For each $r$ choose a family of primitive idempotents
for the modular descent algebra $D_{r,F}$ as explained in \S3.
These are indexed by the $p$-equivalence classes
in ${\rm Part}(r)$.  For any divisor $k$ of $r$ we write $[r,k]$
for the partition of $r$ consisting of $r/k$ parts all equal to
$k$, that is,
$$[r,k] = (k^{r/k}) = (k,k,\dots,k).$$
Furthermore, $e_{[r,k]}$ denotes the primitive idempotent of
$D_{r,F}$ indexed by ${\rm cl}([r,k])$.

We now consider ${\rm cl}([q,q'])$. This
consists of all partitions of $q$ of the form
$$(p^mq',\dots,p^mq',\dots,p^iq',\dots,p^iq',\dots,q',\dots,q'),$$
where the number of occurrences of any of the $p^iq'$ can be $0$.
In particular, $(q) \in {\rm cl}([q,q'])$. Note that
\begin{equation}
{\cal F}_{(q)} = {\cal A}_q = \bigcup_{d \in {\bf D}}
{\cal B}(q/d,d) \cup {\cal A}_q'.
\end{equation}

Let
$${\cal G} = \biggl( \bigcup_{\mu \in {\rm cl}([q,q'])}
{\cal F}_\mu \biggr) \smallsetminus {\cal A}_q'.$$
Thus, by (4.3) and the definition of $C$,
\begin{equation}
\langle {\cal G} \cap L^q(V) \rangle = \bigoplus_{d \in {\bf D}}
L^{q/d}(B_d) = C,
\end{equation}
where $\langle {\cal G} \cap L^q(V) \rangle$
denotes the subspace of $T(V)$ spanned
by ${\cal G} \cap L^q(V)$. If $f \in {\cal G}$ and
$f = a_1 \cdots a_l$, as in the definition of ${\cal F}_q$, then
each factor has degree $p^iq'$ for some $i \in \{0,\dots,m\}$ and
the factors do not belong to ${\cal A}_q'$. Thus, by (4.2),
$$a_1,\dots,a_l \in \bigcup_{d \in {\bf D}} {\cal B}(d).$$
In view of the conditions on the ordering of $\cal A$ we have
$f = \prod_{d \in {\bf D}} f_d$ where the elements $d$ of $\bf D$ are
taken in increasing order and where, for each $d$, $f_d$ is a
product of elements of ${\cal B}(d)$ in increasing order.

Let $d \in {\bf D}$. By (4.1), each element of ${\cal B}(d)$ has
degree divisible by $(q'/d')d$. Therefore every product of elements of
${\cal B}(d)$ has degree $cd$ for some non-negative integer
$c$ divisible by $q'/d'$. For each such $c$, write ${\cal G}(c,d)$
for the set of all elements of $T(V)$ of the form $a_1a_2 \cdots a_l$
with $a_1,\dots,a_l \in {\cal B}(d)$, $a_1 \preccurlyeq \cdots
\preccurlyeq a_l$
and $\deg a_1\cdots a_l = cd$.

Let $f \in {\cal G}$ and write $f = \prod_d f_d$ as above. For each
$d$ write $\deg f_d = c_dd$. Thus $c_d$ is divisible by
$q'/d'$ and $f_d \in {\cal G}(c_d,d)$. Since $\deg f = q$
we have
\begin{equation}
\sum_d c_dd = q.
\end{equation}

Let $S$ be the set of all families $\{c_d: d \in {\bf D}\}$ of
non-negative integers $c_d$ satisfying (4.5) such that $c_d$ is
divisible by $q'/d'$ for all $d$. Note that (4.5) gives
\begin{equation}
c_d \leq q/d = p^{m-j(d)}q'/d'.
\end{equation}
For each $\sigma \in S$, where $\sigma = \{c_d : d \in {\bf D} \}$,
let ${\cal G}(\sigma)$ be the subset of $T(V)$ defined by
\begin{equation}
{\cal G}(\sigma) = \prod_d {\cal G}(c_d,d).
\end{equation}
Thus ${\cal G} = \bigcup_{\sigma \in S} {\cal G}(\sigma)$.

We next consider ${\cal G}(c_d,d)$ for fixed $d$ and $c_d$.
Since $B_d \subseteq T^d(V)$ we may regard $T(B_d)$ as a subalgebra of
$T(V)$, as explained early in \S2. Thus $T^{c_d}(B_d) \subseteq
T^{c_dd}(V)$.  We shall modify ${\cal G}(c_d,d)$ to get a
new set ${\cal G}^*(c_d,d)$.

Consider first the case where $c_d \neq 0$. We form a PBW-basis of
$T^{c_d}(B_d)$ as follows. We start with bases of the Lie powers of
$B_d$ using ${\cal B}(p^iq'/d,d)$ as the basis of $L^{p^iq'/d}(B_d)$
for all $i \in \{j(d),\dots,m\}$. Hence we get a basis of $L(B_d)$
containing ${\cal B}(d)$ which we order in any way subject to
agreement with the order of ${\cal B}(d)$ that we already have.
The ordered products of degree $c_dd$ (as elements of $T(V)$) give
the required PBW-basis of $T^{c_d}(B_d)$. The PBW-filtration
of $T^{c_d}(B_d)$ has factors indexed by the partitions of $c_d$.

Consider the $p$-equivalence class ${\rm cl}([c_d,q'/d'])$. By (4.6),
this class consists of all partitions of $c_d$ of the form
$$(p^{m-j(d)}q'/d',\dots,p^{m-j(d)}q'/d',\dots,
p^iq'/d',\dots,p^iq'/d',\dots,q'/d',\dots,q'/d').$$
The PBW-basis elements corresponding to these partitions are
precisely the elements of ${\cal G}(c_d,d)$.

Since $B_d$ is a module for ${\rm GL}(n,F)$,
the action of ${\rm Sym}(c_d)$ by place
permutations on $T^{c_d}(B_d)$ commutes with the action of
${\rm GL}(n,F)$. Thus $e_{[c_d,q'/d']}T^{c_d}(B_d)$ is a direct
summand of $T^{c_d}(B_d)$ as an $F{\rm GL}(n,F)$-module. By hypothesis, 
$B_d$ is a direct summand of $T^d(V)$, and so  
$e_{[c_d,q'/d']}T^{c_d}(B_d)$ is a direct summand of
$T^{c_dd}(V)$.

For each $f_d \in {\cal G}(c_d,d)$ define $f_d^* =
e_{[c_d,q'/d']}f_d$ and write
\begin{equation}
{\cal G}^*(c_d,d) = \{f_d^* : f_d \in {\cal G}(c_d,d)\}.
\end{equation}
By Proposition 3.5 (applied to the module $B_d$ rather than $V$),
$$\langle {\cal G}^*(c_d,d)\rangle = e_{[c_d,q'/d']}T^{c_d}(B_d).$$
Thus $\langle {\cal G}^*(c_d,d)\rangle $ is an $F{\rm GL}(n,F)$-module
which is a direct summand of $T^{c_dd}(V)$.

If $c_d = 0$ we have ${\cal G}(c_d,d) = \{1\}$. We
define $1^* = 1$. Thus we take ${\cal G}^*(c_d,d)
= \{1\}$ and $\langle {\cal G}^*(c_d,d)\rangle = F.1 = T^0(V)$.

For $f \in {\cal G}(\sigma)$, where $\sigma = \{c_d : d \in {\bf D}\}$
and $f = \prod_d f_d$ as above, we define $f^* = \prod_d f_d^*$.
For each $d$, (3.4) gives that $f_d^*$ is equal to $f_d$ modulo terms
below $f_d$ in the PBW-filtration of $T^{c_d}(B_d)$. It follows that
$f^*$ is equal to $f$ modulo terms below $f$ in the PBW-filtration
of $T^q(V)$. Hence $f$ may be replaced by $f^*$ in the basis
${\cal F}_q$ of $T^q(V)$.

This procedure applied to all $\sigma$ gives $f^*$ for all
$f \in {\cal G} = \bigcup_\sigma {\cal G}(\sigma)$. For $f \in
{\cal F}_q \smallsetminus {\cal G}$ we define $f^* = f$, and write
${\cal F}^*_q = \{f^* : f \in {\cal F}_q\}$. Thus ${\cal F}^*_q$ is
a basis of $T^q(V)$. For each $\lambda \in {\rm Part}(q)$ we
also write ${\cal F}^*_\lambda = \{f^* : f \in {\cal F}_\lambda\}$.
Note that if $f \in {\cal F}_{(q)}$ then $f$ belongs to the bottom
term of the PBW-filtration of $T^q(V)$ and so $f^* = f$.

Recall that ${\cal G}^*(c_d,d)$ is given by (4.8). Similarly, write
$${\cal G}^*(\sigma) = \{f^*:f \in {\cal G}(\sigma)\},$$
for all $\sigma \in S$, and
\begin{equation}
{\cal G}^* = \{f^* : f \in {\cal G}\} = \bigcup_\sigma {\cal G}^*(\sigma).
\end{equation}
By (4.7), if $\sigma = \{c_d: d \in {\bf D}\}$, we have
$${\cal G}^*(\sigma) = \prod_d{\cal G}^*(c_d,d).$$
Thus
\begin{equation}
\langle {\cal G}^*(\sigma)\rangle = \prod_d \langle {\cal G}^*
(c_d,d)\rangle,
\end{equation}
where the right-hand side denotes the span of the product of the sets
$\langle {\cal G}^*(c_d,d)\rangle$ within $T(V)$. Since
${\cal G}^*(\sigma) \subseteq {\cal F}^*_q$, it follows that
${\cal G}^*(\sigma)$ is
a basis of $\langle {\cal G}^*(\sigma)\rangle$. For a similar
reason the spaces $\langle {\cal G}^*(\sigma)\rangle$, $\sigma \in S$,
span their direct sum in $T^q(V)$. Let $U$ be this direct sum,
\begin{equation}
U = \bigoplus_{\sigma \in S} \langle {\cal G}^*(\sigma)\rangle.
\end{equation}
Thus, by (4.9), $U = \langle {\cal G}^* \rangle$. Also, since
${\cal F}^*_{(q)} = {\cal F}_{(q)}$,
$$L^q(V) = \langle ({\cal G} \cap L^q(V)) \cup {\cal A}_q'\rangle
= \langle ({\cal G}^* \cap L^q(V)) \cup {\cal A}_q'\rangle.$$
Hence $U + L^q(V) = \langle {\cal G}^* \cup {\cal A}_q'\rangle$
and $U \cap L^q(V) = \langle {\cal G} \cap L^q(V)\rangle$.
Thus, by (4.4),
\begin{equation}
U \cap L^q(V) = C.
\end{equation}

Since each $\langle {\cal G}^*(c_d,d)\rangle$ is an
$F{\rm GL}(n,F)$-submodule of $T(V)$, so is
$\langle {\cal G}^*(\sigma)\rangle$, by (4.10). Indeed, since
$\langle {\cal G}^*(c_d,d)\rangle$ is a direct summand of
$T^{c_dd}(V)$ for all $d$, it follows that
$\langle {\cal G}^*(\sigma)\rangle$ is a direct summand of $T^q(V)$.
Thus, by Lemma 2.9, each $\langle {\cal G}^*(\sigma)\rangle$ is
injective as an $S_F(n,q)$-module. Thus, by (4.11), $U$ is
injective.

Let $\phi:T^q(V) \to e_{[q,q']}T^q(V)$ be the mapping
given by $\phi(w) = e_{[q,q']}w$ for all $w \in T^q(V)$.
Thus (see \S3) $\phi$ is a homomorphism of $F{\rm GL}(n,F)$-modules;
hence a homomorphism of $S_F(n,q)$-modules.
As shown above, $U+L^q(V)$ is the subspace of $T^q(V)$ spanned by
${\cal G}^* \cup {\cal A}_q'$. However,
$${\cal G}^* \cup {\cal A}_q' = \bigcup_{\mu \in {\rm cl}([q,q'])}
{\cal F}_\mu^*.$$
Hence, by Proposition 3.5,
the restriction of $\phi$ is a vector space isomorphism from
$U + L^q(V)$ to $e_{[q,q']}T^q(V)$. Thus it is an
$S_F(n,q)$-module isomorphism.
Therefore, by Lemma 2.9, $U+L^q(V)$ is injective and projective.
(Hence it is a direct summand of $T^q(V)$.) However, by (4.12),
$$L^q(V)/C \cong (U+L^q(V))/U.$$
Since $U$ is injective, it follows that $L^q(V)/C$ is injective
and projective. This completes the proof of Proposition 4.1.
\hfill{$\square$}

\vspace{0.2in}
We can now state our key result. We refer to \S2
for the notation $V^{(n)}$ and definitions concerning uniform families
of modules.

\vspace{0.2in}
{\sc Theorem} 4.2.  {\it Let\/ $F$ be a field and\/ $k$ a positive
integer not divisible by the characteristic of\/ $F$. For each
positive integer\/ $s$ there exists a uniform family\/
$\{B_{sk}^{(n)} : n \in {\mathbb N}\}$ of modules of degree\/ $sk$
such that, for all\/ $n$, $B_{sk}^{(n)}$ is a submodule of\/
$L^{sk}(V^{(n)})$ and a direct summand of\/ $T^{sk}(V^{(n)})$ and
$$L^k(V^{(n)}) \oplus L^{2k}(V^{(n)}) \oplus L^{3k}(V^{(n)}) \oplus \cdots
= L(B_k^{(n)}) \oplus L(B_{2k}^{(n)}) \oplus L(B_{3k}^{(n)}) \oplus
\cdots \; .$$
}

\vspace{0.01in}
Before proving Theorem 4.2 we show that it implies our main results
for an arbitrary group $G$.

\vspace{0.2in}
{\sc Theorem} 4.3. {\it Let\/ $F$ be a field and\/ $k$ a positive
integer not divisible by the characteristic of\/ $F$. Let\/ $G$ be
a group and\/ $V$ a finite-dimensional\/ $FG$-module. For each
positive integer\/ $s$ there is a submodule\/ $B_{sk}$ of\/
$L^{sk}(V)$ such that\/ $B_{sk}$ is a direct summand of\/ $T^{sk}(V)$
and
$$L^{sk}(V) = \bigoplus_{c \mid s} L^{s/c}(B_{ck}).$$}

\vspace{0.01in}
{\it Proof}. Let $n = \dim V$. By choice of a basis of $V$, the
representation of $G$ on $V$ gives a homomorphism $\theta :
G \to {\rm GL}(n,F)$. We can regard $V$ as the natural module for
${\rm GL}(n,F)$. For each positive integer $s$ define
$B_{sk} = B_{sk}^{(n)}$ where $B_{sk}^{(n)}$ is as given by
Theorem 4.2. The action of $G$ on $T(V)$ is the composite of
$\theta$ and the action of ${\rm GL}(n,F)$ on $T(V)$. Hence,
for each $s$, $B_{sk}$ is an $FG$-submodule of $L^{sk}(V)$ and
a direct summand of $T^{sk}(V)$, and $\bigoplus_s L^{sk}(V) =
\bigoplus_s L(B_{sk})$. By comparing terms of degree $sk$ we obtain
$L^{sk}(V) = \bigoplus_{c \mid s} L^{s/c}(B_{ck})$.
\hfill{$\square$}

\vspace{0.2in}
{\sc Theorem 4.4 (Decomposition Theorem)}.
{\it Let\/ $F$ be a field of prime
characteristic\/ $p$.  Let\/ $G$ be a group and\/ $V$ a
finite-dimensional\/ $FG$-module.  For each positive integer\/ $r$
there is a submodule\/ $B_r$ of\/ $L^r(V)$ such that\/ $B_r$ is
a direct summand of\/ $T^r(V)$ and,
for\/ $k$ not divisible by\/ $p$ and\/ $m \geq 0$,
$$L^{p^mk}(V) = L^{p^m}(B_k) \oplus L^{p^{m-1}}(B_{pk}) \oplus
\cdots \oplus L^1(B_{p^mk}).$$}

\vspace{0.01in}
{\it Proof}. For each positive integer $k$ not divisible by $p$,
take $B_k$, $B_{pk}$, $B_{p^2k}$, \dots\ as given by Theorem 4.3.
The result follows. \hfill{$\square$}

\vspace{0.2in}
{\it Proof of Theorem\/} 4.2.
We shall consider various families $\{U^{(n)}:n \in {\mathbb N}\}$
where, for each $n$, $U^{(n)}$ is a subspace of $T(V^{(n)})$.
When $n$ is regarded as fixed
but arbitrary we generally simplify the
notation by writing $U$ instead of $U^{(n)}$. In particular, we often
write $V$ for $V^{(n)}$.

Let $Q = Q^{(n)} = L^k(V) \oplus L^{2k}(V) \oplus L^{3k}(V) \oplus
\cdots\;$.
Clearly $Q$ is a subalgebra of $L(V)$ and
$Q = \bigoplus_{r \geq 1} Q_r$ where $Q_r = Q \cap L^r(V)$ for
all $r$. Of course, $Q_r = 0$ unless $r$ is divisible by $k$ and
$Q_{sk} = L^{sk}(V)$ for all $s$.
For each $r \geq 1$, let $Q(r)$ denote the
subalgebra of $Q$ generated by $Q_1 \oplus \cdots \oplus Q_r$. Also,
write $Q(0) = 0$.

For each $r \geq 1$, let $W_r$ (or $W_r^{(n)}$) be any subspace of
$Q_r$ with the property that
$Q_r = (Q(r-1) \cap Q_r) \oplus W_r$. The choice of these subspaces
is arbitrary at present: later we make a more careful choice.
By Lemma 2.1, we have $Q = L(W_1 \oplus W_2 \oplus \cdots)$.
Clearly $W_r = 0$ unless $r = sk$ for some $s$.
Hence $Q = L(W_k \oplus W_{2k} \oplus W_{3k} \oplus \cdots )$.
Also, $W_k = L^k(V)$.

For each positive integer $s$ we construct a subspace
$B_{sk}$ of $L^{sk}(V)$ such that
$Q = L(B_k) \oplus L(B_{2k}) \oplus \cdots \;$.
For each $s$, $B_{sk}$ depends only on $W_k, W_{2k}, \dots, W_{sk}$.

Step 1. Define $B_k = W_k$ and $C_k = W_{2k} \oplus W_{3k} \oplus
\cdots \;$. By Lemma 2.2,
\begin{equation}
Q = L(B_k \oplus C_k) = L(B_k) \oplus L(C_k \wr B_k).
\end{equation}
Also, by Corollary 2.3, $C_k \wr B_k$ is the direct sum of the subspaces
in ${\cal E}_1$, where
$${\cal E}_1 = \{[W_{ik}, \underbrace{B_k, \dots, B_k}_r] :
i>1,\; r \geq 0\}.$$

Step 2. Let $B_{2k}$ be the sum of all elements of ${\cal E}_1$
of degree $2k$ (so that, in fact, $B_{2k} = W_{2k}$)
and let $C_{2k}$ be the sum of all
elements of ${\cal E}_1$ of degree greater than $2k$.
By Lemma 2.2,
$$L(C_k \wr B_k) = L(B_{2k}) \oplus L(C_{2k} \wr B_{2k}).$$
Also, by Corollary 2.3,
$C_{2k}\wr B_{2k}$ is the direct sum of the
subspaces in ${\cal E}_2$, where
$${\cal E}_2 = \{ [E,\underbrace{B_{2k},\dots,B_{2k}}_r] :
E \in {\cal E}_1,\;
\deg E > 2k,\; r \geq 0\}.$$

We continue in this way. 
In Step $s$ (where $s > 1$) we define $B_{sk}$ to be the sum of
all elements of ${\cal E}_{s-1}$ of degree
$sk$ and $C_{sk}$ to be the sum of all elements of ${\cal E}_{s-1}$
of degree greater than $sk$. By Lemma 2.2 and
Corollary 2.3,
\begin{equation}
L(C_{(s-1)k} \wr B_{(s-1)k}) =
L(B_{sk}) \oplus L(C_{sk} \wr B_{sk})
\end{equation}
and $C_{sk} \wr B_{sk}$ is the direct sum of the subspaces in ${\cal E}_s$,
where
$${\cal E}_s = \{[E,\underbrace{B_{sk},\dots,B_{sk}}_r] :
E \in {\cal E}_{s-1},\;
\deg E > sk,\; r \geq 0\}.$$
In particular, every subspace in ${\cal E}_s$ has the form
$$[W_{ik},B_{s(1)k},\dots,B_{s(r)k}]$$
where $i>1$, $r \geq 0$ and $1 \leq s(1) \leq \cdots \leq s(r) \leq s$.

By (4.13), (4.14) and considerations of degree, we obtain
\begin{equation}
Q = L(B_k) \oplus L(B_{2k}) \oplus L(B_{3k}) \oplus \cdots \; .
\end{equation}

This will give the statement of the theorem
if we can show that the families $\{W_{sk}^{(n)}\}$
may be chosen so that the resulting families $\{B_{sk}^{(n)}\}$
are uniform families of modules with $B_{sk}^{(n)}$ a direct
summand of $T^{sk}(V^{(n)})$ for all $n$ and $s$.

Comparing terms of degree $sk$ in (4.15) we have, for all $s \geq 1$,
\begin{equation}
L^{sk}(V) = Q_{sk} = \bigoplus_{c \mid s} L^{s/c}(B_{ck}).
\end{equation}

As at the end of \S2, let ${\cal T}_r^{(n)}$ (or ${\cal T}_r$
when $n$ is understood) denote the class of all finite-dimensional
$S_F(n,r)$-modules that are direct sums of modules each of which is
isomorphic to a direct summand of $T^r(V^{(n)})$. By Lemma 2.10
it suffices to obtain uniform families $\{B_{sk}^{(n)}\}$
such that $B_{sk}^{(sk)} \in {\cal T}_{sk}^{(sk)}$ for all $s$.

We start with arbitrary families $\{W_{sk}^{(n)}\}$ and then
refine the choice of $\{W_k^{(n)}\}$, $\{W_{2k}^{(n)}\}$, \dots,
successively, so that the resulting families $\{W_{sk}^{(n)}\}$ and
$\{B_{sk}^{(n)}\}$ are uniform families of modules satisfying
$W_{sk}^{(n)}, B_{sk}^{(n)} \in {\cal T}_{sk}^{(n)}$ for all
$n$ and $s$.

There is no choice possible for $W_k$: we have $W_k =
L^k(V)$. Step 1 of the construction described above gives
$B_k = W_k$. Thus the families $\{W_k^{(n)}\}$ and
$\{B_k^{(n)}\}$ are uniform families of modules of degree $k$.
By a well-known result (see, for example, [{\bf 13}, \S3.1]),
$L^k(V)$ is a direct summand of $T^k(V)$, because $k$ is not divisible
by the characteristic of $F$.
Thus $W_k \in {\cal T}_k$ and $B_k \in {\cal T}_k$.

Suppose that $s \geq 2$ and that we have chosen the families
$\{W_{ik}^{(n)}\}$ for all $i < s$ so that $\{W_{ik}^{(n)}\}$ and
$\{B_{ik}^{(n)}\}$ are uniform families of modules
belonging to ${\cal T}_{ik}$. We now consider Step $s$.

By definition,
$B_{sk}$ is the sum of all elements of ${\cal E}_{s-1}$ of degree
$sk$. One such term is $W_{sk}$. The other
terms have the form $[W_{ik},B_{s(1)k},\dots,B_{s(r)k}]$
with $r > 0$. Let $U_{sk}$ denote the sum of these
latter terms. Thus $B_{sk} = U_{sk} \oplus W_{sk}$.
Consider a summand
$[W_{ik},B_{s(1)k},\dots,B_{s(r)k}]$ of $U_{sk}$. Since $r > 0$,
we have $i < s$ and $s(1),\dots,s(r) < s$. Hence
$\{W_{ik}^{(n)}\}$, $\{B_{s(1)k}^{(n)}\}$, \dots, $\{B_{s(r)k}^{(n)}\}$
are uniform families of modules belonging to ${\cal T}_{ik}$,
${\cal T}_{s(1)k}$, \dots, ${\cal T}_{s(r)k}$, respectively.

Repeated use of Corollary 2.3 gives a vector space isomorphism
$$[W_{ik},B_{s(1)k}, \dots, B_{s(r)k}] \cong
W_{ik} \otimes B_{s(1)k} \otimes \cdots \otimes B_{s(r)k}.$$
It is easy to check, by Lemma 2.4,
that these spaces are $S_F(n,sk)$-modules and the isomorphism is a
module isomorphism. It follows that
$[W_{ik},B_{s(1)k},\dots,B_{s(r)k}] \in {\cal T}_{sk}$.
Also, by Lemma 2.7, the family $\{[W_{ik}^{(n)},B_{s(1)k}^{(n)},
\dots,B_{s(r)k}^{(n)}]\}$ is uniform. It follows that $\{U_{sk}^{(n)}\}$ is
a uniform family of modules belonging to ${\cal T}_{sk}$.

By (4.16),
\begin{equation}
L^{sk}(V) = Q_{sk} = \bigoplus_{c \mid s,\, c<s} L^{s/c}(B_{ck})
\oplus U_{sk} \oplus W_{sk}.
\end{equation}
Note that
$$\bigoplus_{c \mid s,\, c<s} L^{s/c}(B_{ck}) \oplus U_{sk}
\subseteq Q(sk-1) \cap Q_{sk}$$
and, by the defining property of $W_{sk}$,
\begin{equation}
Q_{sk} = (Q(sk-1) \cap Q_{sk}) \oplus W_{sk}.
\end{equation}
Thus, by (4.17),
\begin{equation}
\bigoplus_{c \mid s,\, c<s} L^{s/c}(B_{ck})
\oplus U_{sk} = Q(sk-1) \cap Q_{sk}.
\end{equation}

We next prove that the module
\begin{equation}
L^{sk}(V^{(sk)})\bigg/ \bigoplus_{c \mid s,\, c<s} L^{s/c}(B_{ck}^{(sk)})
\end{equation}
is injective and projective as an $S_F(sk,sk)$-module. This is clear if
$F$ has characteristic $0$ because in that case $S_F(sk,sk)$ is
semisimple, by [{\bf 16}, (2.6e)]. Thus we may suppose that $F$ has prime
characteristic $p$. By hypothesis, $p \nmid k$. For each positive
integer $r$ not divisible by $k$, define $B_r^{(sk)} = 0$. In
Proposition 4.1, take $n=q=sk$, $V = V^{(sk)}$ and $B_d = B_d^{(sk)}$
for every proper divisor $d$ of $q$. The hypotheses of Proposition
4.1 follow from (4.16). Thus the module (4.20)
is injective and projective, as required.

By (4.17), this module has a submodule
$$\biggl( \bigoplus_{c \mid s,\, c<s}
L^{s/c}(B_{ck}^{(sk)}) \oplus U_{sk}^{(sk)} \biggr)
\bigg/ \bigoplus_{c \mid s,\, c<s} L^{s/c}(B_{ck}^{(sk)})$$
isomorphic to $U_{sk}^{(sk)}$. By Lemma 2.9, $U_{sk}^{(sk)}$
is injective. It follows that
$$L^{sk}(V^{(sk)})\bigg/\biggl(\bigoplus_{c \mid s,\, c<s}
L^{s/c}(B_{ck}^{(sk)}) \oplus U_{sk}^{(sk)}\biggr)$$
is injective and projective. Hence there exists $Z_{sk}^{(sk)}
\in {\cal T}_{sk}^{(sk)}$ such that
$$L^{sk}(V^{(sk)}) = \bigoplus_{c \mid s,\, c<s} L^{s/c}(B_{ck}^{(sk)})
\oplus U_{sk}^{(sk)} \oplus Z_{sk}^{(sk)}.$$
By Lemma 2.5, $Z_{sk}^{(sk)}$ extends to a uniform family
$\{Z_{sk}^{(n)}\}$ of modules of degree $sk$. Also, by Lemma 2.10,
$Z_{sk}^{(n)} \in {\cal T}_{sk}^{(n)}$ for all $n$.
Furthermore, by Lemma 2.7, $\{L^{s/c}(B_{ck}^{(n)})\}$ is uniform,
for all proper divisors $c$ of $s$,
and we have seen that $\{U_{sk}^{(n)}\}$ is uniform.
Hence, by Lemma 2.6, we have (for arbitrary $n$)
$$L^{sk}(V) = \bigoplus_{c \mid s,\, c<s} L^{s/c}(B_{ck})
\oplus U_{sk} \oplus Z_{sk}.$$
Therefore, in view of (4.19),
$$Q_{sk} = L^{sk}(V) = (Q(sk-1) \cap Q_{sk}) \oplus Z_{sk}.$$
However, (4.18) is the defining property of $W_{sk}$. Thus we
may replace our arbitrary choice of $W_{sk}$ by taking
$W_{sk} = Z_{sk}$. Recall that $B_{sk} = U_{sk} \oplus
W_{sk}$. Thus $\{W_{sk}^{(n)}\}$ and $\{B_{sk}^{(n)}\}$ are
uniform families of modules belonging to ${\cal T}_{sk}$.

This completes the proof of Theorem 4.2. \hfill{$\square$}

\nopagebreak
\vspace{0.4in}
{\it Acknowledgements.} The authors are grateful to Dr K. Erdmann and
Dr R. St\"ohr for their useful comments while this work was in progress.

\vspace{0.4in}
\begin{center}
{\it References}
\end{center}

\newcounter{euro}
\begin{list}{\bf \arabic{euro}.\hspace{0.1in}}{\parsep=0pt
\itemsep=0pt \usecounter{euro}}

\item {\sc D. Blessenohl} and {\sc H. Laue}, `The module structure of
Solomon's descent algebra', {\it J. Austral.\ Math.\ Soc.} 72 (2002) 
317--333.
\item {\sc N. Bourbaki}, {\it Lie groups and Lie algebras}, Part I: Chapters
1--3 (Hermann, Paris, 1975).
\item {\sc A. J. Brandt}, `The free Lie ring and Lie representations of the
full linear group', {\it Trans.\ Amer.\ Math.\ Soc.} 56 (1944) 
528--536.
\item {\sc R. M. Bryant}, `Free Lie algebras and Adams operations',
{\it J. London Math.\ Soc.} (2) 68 (2003) 355--370.
\item {\sc R. M. Bryant}, `Modular Lie representations of groups of
prime order', {\it Math.\ Z.} 246 (2004) 603--617.
\item {\sc R. M. Bryant}, `Modular Lie representations of finite groups',
{\it J. Austral.\ Math.\ Soc.} 77 (2004) 401--423.
\item {\sc R. M. Bryant, L. G. Kov\'acs} and {\sc R. St\"ohr},
`Lie powers of modules for groups of prime order', {\it Proc.\ London
Math.\ Soc.} (3) 84 (2002) 343--374.
\item {\sc R. M. Bryant} and {\sc I. C. Michos}, `Lie powers of free modules
for certain groups of prime power order', {\it J. Austral.\ Math.\ 
Soc.} 71 (2001) 149--158.
\item {\sc R. M. Bryant} and {\sc M. Schocker}, `Lie powers and Lie
resolvents', in preparation.
\item {\sc R. M. Bryant} and {\sc R. St\"{o}hr}, `Lie powers in prime
degree', {\it Quart.\ J. Math.}, to appear.
\item {\sc H. Cartan} and {\sc S. Eilenberg}, {\it Homological algebra\/} 
(Princeton University Press, Princeton, 1956).
\item {\sc S. Donkin}, `On tilting modules for algebraic groups',
{\it Math.\ Z.} 212 (1993) 39--60.
\item {\sc S. Donkin} and {\sc K. Erdmann}, `Tilting modules, symmetric
functions, and the module structure of the free Lie algebra',
{\it J. Algebra\/} 203 (1998) 69--90.
\item {\sc K. Erdmann} and {\sc M. Schocker}, `Modular Lie powers and
the Solomon descent algebra', preprint.
\item {\sc A. M. Garsia} and {\sc C. Reutenauer}, `A decomposition of
Solomon's descent algebra', {\it Adv.\ Math.} 77 (1989) 
189--262.
\item {\sc J. A. Green}, {\it Polynomial representations of\/ ${\rm GL}_n$},
Lecture Notes in Mathematics 830 (Springer, Berlin, 1980).
\item {\sc B. Hartley} and {\sc R. St\"ohr}, `Homology of higher relation
modules and torsion in free central extensions of groups',
{\it Proc.\ London Math.\ Soc.} (3) 62 (1991) 325--352.
\item {\sc N. Jacobson}, {\it Lie algebras\/} (Interscience, New York, 1962).
\item {\sc A. A. Klyachko}, `Lie elements in the tensor algebra',
{\it Siberian Math.\ J.} 15 (1974) 914--921 (1975).
\item {\sc W. Kraskiewicz} and {\sc J. Weyman}, `Algebra of coinvariants
and the action of a Coxeter element', {\it Bayreuth.\ Math.\ Schr.}
63 (2001) 265--284. 
\item {\sc C. Reutenauer}, {\it Free Lie algebras\/} (Clarendon Press,
Oxford, 1993).

\end{list}

\vspace{0.1in}
\baselineskip=16pt
\begin{tabbing}
{\it R. M. Bryant  \hspace{2.0in}} \=
{\it M. Schocker}\\
{\it School of Mathematics} \> {\it Department of Mathematics}\\
{\it University of Manchester} \> {\it University of Wales Swansea}\\
{\it PO Box 88} \> {\it Singleton Park}\\
{\it Manchester M60 1QD} \> {\it Swansea SA2 8PP}\\
roger.bryant@manchester.ac.uk \> m.schocker@swansea.ac.uk
\end{tabbing}

\end{document}